\documentclass[11pt,a4paper,oneside]{article}
\usepackage[displaymath, mathlines]{lineno}
\usepackage{authblk}
\usepackage{xr}
\usepackage{natbib}
\usepackage{amsmath}
\usepackage{mathrsfs}  
\usepackage{amsfonts}
\usepackage{amssymb}
\usepackage{amsthm}
\usepackage{graphicx}
\usepackage{subcaption}
\usepackage{soul, color, colortbl}
\usepackage{booktabs}
\usepackage[margin=1.25in]{geometry}
\usepackage{enumerate}
\usepackage[shortlabels]{enumitem}
\usepackage{url}
\usepackage{afterpage}
\usepackage{multirow}

\newtheorem{assumption}{Assumption}
\newtheorem{theorem}{Theorem}
\newtheorem{corollary}[theorem]{Corollary}

\newtheorem{remark}{Remark}
\numberwithin{assumption}{section}
\numberwithin{theorem}{section}
\numberwithin{proposition}{section}
\numberwithin{remark}{section}
\numberwithin{lemma}{section}

\newcommand{\bref}[1]{(\ref{#1})}

\newcolumntype{L}[1]{>{\raggedright\let\newline\\\arraybackslash\hspace{0pt}}m{#1}}
\newcolumntype{C}[1]{>{\centering\let\newline\\\arraybackslash\hspace{0pt}}m{#1}}
\newcolumntype{R}[1]{>{\raggedleft\let\newline\\\arraybackslash\hspace{0pt}}m{#1}}

\begin{document}

\title{ {\bf Robust regression estimation and inference in the presence of cellwise and casewise contamination}}
\author{{\bf Andy Leung}}
\author{{\bf Hongyang Zhang}}
\author{{\bf Ruben H. Zamar}}
\affil{{\small Department of Statistics, University of British Columbia, 3182-2207 Main Mall, Vancouver, British Columbia V6T 1Z4, Canada}}
\renewcommand\Authands{, }

\maketitle

\begin{abstract}
Cellwise outliers are likely to occur together with casewise outliers in modern datasets of relatively large dimension. Recent work has shown that traditional robust regression methods may fail when applied to such datasets.  We propose a new robust regression procedure to deal with casewise and cellwise outliers. The proposed method, called three-step regression, proceeds as follows:  first, it uses a consistent univariate  filter, that is, a procedure that  flags and eliminates extreme cellwise outliers; second, it  applies  a robust estimator of multivariate location and scatter to the filtered data to down-weight  casewise outliers; third, it computes robust regression coefficients from the estimates obtained in the second step. The three-step estimator is consistent and asymptotically normal at the central model under some assumptions on the tails of the  distributions of the continuous covariates. The estimator is extended to handle both continuous and dummy covariates using an iterative algorithm. Extensive simulation results show that the three-step estimator is resilient to cellwise outliers. It also performs  well under casewise contamination when compared to traditional high breakdown point estimators. 
\end{abstract}

\section{Introduction}

The vast majority of procedures for robust  linear regression are based on the classical Tukey--Huber contamination model (THCM) in which a relatively small fraction of cases may be  contaminated. High breakdown point affine equivariant estimators such as least trimmed squares  \citep{rousseeuw:1984b}, S-regression \citep{rousseeuw:1984a} and MM-regression \citep{yohai:1985} proceed by down-weighting   outlying cases, which makes sense and works well in practice, under THCM. However, in some applications, the contamination mechanism may be different in that random cells in a data table (with rows as cases and columns as variables)  are independently contaminated. In this paradigm, a small fraction of random cellwise outliers could  propagate to a relatively large fraction of cases, breaking down  classical high breakdown point affine equivariant estimators  \citep[see][]{alqallaf:2009}. Since cellwise and casewise outliers may co-exist  in some applications,  our goal in this paper is to develop a method for robust regression estimation and inference that can deal with both cellwise and casewise outliers.

There is a vast literature on robust regression for casewise outliers, but only a scant literature for  cellwise outliers and none for both types of outliers in the regression context.  Recently, \citet{ollerer:2015} combined the ideas of  coordinate descent algorithm \citep[called the shooting algorithm in][]{fu:1998} and simple S-regression \citep{rousseeuw:1984a} to propose an estimator  called the shooting S. The shooting S-estimator  assigns individual weight to each cell in the data table to  handle cellwise outliers in the regression context. The shooting S-estimator is  robust against cellwise outliers and vertical response outliers. 

In this paper, we propose a three-step regression estimator  which combines the ideas of  filtering cellwise outliers and robust regression via covariance matrix estimate \citep{maronna:1986, croux:2003}, namely 3S-regression estimator. By filtering, here we mean  detecting outliers and replacing them by missing values as in \citet{agostinelli:2014}.
Our estimator proceeds as follows:  first, it uses a univariate  filter to detect and eliminate extreme cellwise outliers in order to control the effect of  outliers propagation; second, it  applies  a robust estimator of multivariate location and scatter to the filtered data to down-weight  casewise outliers; third, it computes robust regression coefficients from the estimates obtained in the second step. With the choice of a filter that has simultaneous good sensitivity (is capable of filtering outliers) and good specificity (can preserve all or most of the clean data), the resulting estimator can be resilient to both cellwise and casewise outliers; furthermore, it attains consistency and asymptotic normality for clean data. In this regards, we propose a filter that is consistent under some assumptions on the tails of the covariates distributions. By consistent filter, we mean a filter that asymptotically can preserve all the data when they are clean. 

The rest of the paper is organized as follows. In Section \ref{sec:filter}, we introduce a family of consistent filters. In Section \ref{sec:3S-regress}, we introduce 3S-regression.  In Section \ref{sec:asymptotics}, we show some asymptotic properties of 3S-regression. In Section \ref{sec:simulation}, we evaluate the performance of 3S-regression in an extensive simulation study. In Section \ref{sec:realdata}, we analyze a real data set with cellwise and casewise outliers. In Section \ref{sec:conclusion}, we conclude with some remarks. We also provide a document referred to as ``supplementary material'', which  contains all the proofs, additional simulation results, and other related material.

\section{Consistent filter}\label{sec:filter}

Filtering is a method for pre-processing  data in order to control the effect of potential cellwise outliers. In this paper, we pre-process the data by flagging outliers and replacing them by missing values, NAs. This method of filtering has recently been  used for robust estimation of multivariate location and scatter \citep{danilov:2010,agostinelli:2014} and for clustering \citep{farcomeni:2014a, farcomeni:2014b}.  Also, \citet{farcomeni:2015} proposed a procedure to determine a data-driven choice for  the number of filtered cells to increase the efficiency of the  estimator.

Consistent filters are ones that do not filter good data points asymptotically.  \citet{gervini:2002} introduced a consistent filter for  normal residuals in regression estimation to achieve a  fully-efficient robust regression estimator. 
Consistent filters are desirable because their good asymptotic properties are shared by the following-up estimation procedure.  In this paper, we introduce a new family of consistent filters for univariate data.

Consider a random variable $X$ with a continuous distribution function $G(x)$. We define the scaled upper and lower tail distributions of  $G(x)$ as follows:
\begin{linenomath} 
\begin{equation}\label{eq:tail}
\begin{aligned}
F^u(t) &= P_{G}\left(  \frac{X - \eta^u}{\text{med}( X - \eta^u | X > \eta^u)} \le t| X > \eta^u \right) \quad \text{and} \\ 
F^l(t) &= P_{G}\left(  \frac{\eta^l - X}{\text{med}(\eta^l - X | X < \eta^l)} \le t| X < \eta^l \right).
\end{aligned}
\end{equation}
\end{linenomath} 
Here, $\text{med}$ stands for median, $\eta^u = G^{-1}(1 - \alpha)$, $\eta^l = G^{-1}(\alpha)$, and $0 < \alpha < 0.5$. We use $\alpha = 0.20$, but other choices could be considered.  To simplify the notation, we  set $s^u = \text{med}(X - \eta^u | X > \eta^u)$ and $s^l = \text{med}(\eta^l - X | X < \eta^l)$. Alternatively, a combined tails approach could be used for symmetric distributions as in \citet{gervini:2002}. 

Let $\{X_1, \dots, X_n\}$ be a random sample from $G$, and let  $X_{(1)} < X_{(2)} < \dots < X_{(n)}$ be the corresponding order statistics.  Consistent estimators for $(\eta^u,  s^u, \eta^l,  s^l)$ are given by
\begin{linenomath} 
\[
\begin{aligned}
\hat \eta^u_{n} &= \hat G^{-1}_n(1- \alpha), \quad \hat s^u_n = \text{med}(\{ X_i - \hat\eta^u_n | X_i > \hat\eta^u_n \} ),\\
\hat \eta^l_{n} &= \hat G^{-1}_n(\alpha), \quad \hat s^l_n = \text{med}(\{ \hat\eta^l_n - X_i | X_i < \hat\eta^l_n \} ),
\end{aligned}
\]
\end{linenomath} 
where $\hat G_n^{-1}(a) = X_{(\lceil na \rceil)}$, $0 < a < 1$, is the empirical quantile and $\text{med}(\{Y_1, \dots, Y_m\}) =  Y_{(\lceil m/2 \rceil)}$ is the sample median  (see Lemma 1.1 in the supplementary material for a proof of the consistency for $\hat s^u_n$ and $\hat s^l_n$). The  empirical distribution functions for the scaled upper and lower tails in \bref{eq:tail} are now given by
\begin{linenomath} 
 \[
\begin{aligned}
\hat F^u_{n}(t) &= \frac{\sum_{i=1}^n I(0 < (X_{i} - \hat\eta^u_{n})/\hat s^u_n \le t)}{ \sum_{i=1}^n I( X_i > \hat\eta^u_{n})} \quad \text{and}\\
\hat F^l_{n}(t) &= \frac{\sum_{i=1}^n I(0 < ( \hat\eta^l_{n} -  X_{i})/\hat s^l_n \le t)}{ \sum_{i=1}^n I( X_{i} < \hat\eta^l_{n})}.
\end{aligned}
\] 
\end{linenomath} 

Upper and lower tails outliers can be flagged by comparing the empirical distribution functions for the scaled tails with their expected distributions. We assume that aside from  contamination, $F^u$ and $F^l$ decay exponentially fast or faster. Let $\{a\}^+ = \max(0, a)$ denote the positive part of $a$. 
Then, we define the proportions of flagged upper and lower tails outliers  by 
\begin{linenomath} 
\[
\begin{aligned}
\hat d^u_{n} &= \sup_{t \ge t_0} \left\{ F_0(t) - \hat F^u_{n}(t) \right\}^+ \quad \text{and} \quad 
\hat d^l_{n}&= \sup_{t \ge t_0} \left\{ F_0(t) - \hat F^l_{n}(t) \right\}^+,
\end{aligned}
\]
\end{linenomath} 
where $F_0(t) = 1 - \exp(-\log(2) t)$ and $t_0= 1/\log(2)$. 
When $X - \eta^u | X> \eta^u$ is exponentially distributed with a rate of $\lambda^u > 0$, the standardized tail $(X - \eta^u)/s^u | X> \eta^u$ have exponential distribution with a rate of $\log(2)$,  leading to our choice of $F_0(t)$ and $t_0$.  
 Finally, we filter  $\hat d^u_{n} \times 100\%$ of the most extreme points in the upper tail $\{ X_{i} |  X_{i} > \hat\eta^u_{n} \}$, and filter  $\hat d^l_{n} \times 100\%$ of the most extreme points in the lower tail $\{ X_{i} | X_{i} < \hat\eta^l_{n} \}$. Equivalently, setting
\begin{linenomath} 
\[
\begin{aligned}
\hat t^u_{n} &= \min \left\{t:  \hat F^u_{n}(t) \ge 1 - \hat d^u_{n} \right\} \quad
\text{and} \quad 
\hat t^l_{n}&= \min \left\{ t: \hat F^l_{n}(t) \ge 1 - \hat d^l_{n} \right\},
\end{aligned}
\]
\end{linenomath} 
we filter $X_i$'s with $X_i < \hat\eta^l_n - \hat s^l_n\hat t^l_n$ or $X_i > \hat\eta^u_n + \hat s^u_n\hat t^u_n$.

We tried several heavy tail models for $F_0(t)$ including Pareto distributions with different tail indexes, and we found that the chosen exponential model strikes a good balance between the robustness and consistency  of the filtering procedure. 

Theorem \ref{thm:filter} (proved in the supplementary material)  below shows that our filter is consistent under the following assumption on the tails of $G(x)$.

\begin{assumption}\label{assum:tails}
$G(x)$ is continuous, and 
$F^u(t)$ and $F^l(t)$ satisfy the following:
\begin{linenomath} 
\[
F_0(t) - F^u(t) \le 0, \quad t \ge t_0
 \quad \text{and} \quad 
F_0(t) - F^l(t) \le 0, \quad t \ge t_0.
\]
\end{linenomath} 
\end{assumption}

\begin{theorem}\label{thm:filter}
Suppose that Assumption \ref{assum:tails} holds for $G(x)$. 
Then, $\hat d^u_n \to 0$ a.s. and $\hat d^l_n \to 0$ a.s.
\end{theorem}

In practice, the  distributions $F^u(t)$ and $F^l(t)$ are  unknown. To allow for some flexibility, Assumption 2.1 does not completely specify $F^u(t)$ and $F^l(t)$, but it only requires that their upper tails are as heavy as  or lighter than   the upper  tail  of $F_0(t)$.

\section{Three-step regression}\label{sec:3S-regress}

\subsection{The estimator}

Consider the model
\begin{linenomath} 
\begin{equation}\label{eq:linear-model}
	Y_i = \alpha + \pmb X_i^t\pmb\beta + \varepsilon_i
\end{equation}
\end{linenomath} 
for $i=1,\dots,n$, where the error terms $\varepsilon_i$ are i.i.d. and independent of the   covariates $\pmb X_i = (X_{i1},\dots,X_{ip})^t$. 
The least squares (LS) estimates $(\hat{\alpha}_{LS}, \hat{\pmb\beta}^t_{LS})$ are defined as the minimizers of the sum squares of residuals,
\begin{linenomath} 
\[
   (\hat{\alpha}_{LS}, \hat{\pmb\beta}^t_{LS} ) = \underset{(\alpha, \pmb\beta^t) \in \mathbb R^{(p+1)}}{\arg \min} \sum_{i=1}^n ( Y_i - \alpha - \pmb X_i^t\pmb\beta )^2.
\]
\end{linenomath} 
The solution to this problem is explicit:
\begin{linenomath} 
\begin{equation}\label{eq:sol-LS}
\begin{aligned}
\hat{\pmb\beta}_{LS} &= \hat{\pmb\Sigma}_{xx}^{-1} \hat{\pmb\Sigma}_{xy}, \\
 \hat{\alpha}_{LS} &= \hat{\mu}_y - \hat{\pmb\mu}_x^t \hat{\pmb\beta}_{LS}.
\end{aligned}
\end{equation}
\end{linenomath} 
Here,  $\hat{\pmb\Sigma}_{xx}, \hat{\pmb\Sigma}_{xy}, \hat{\mu}_y$, and $\hat{\pmb\mu}_x$ are the  components of the empirical covariance matrix and mean:
\begin{linenomath} 
\begin{equation}\label{eq:partition}
\hat{\pmb\Sigma} = \left( \begin{array}{cc} \hat{\pmb\Sigma}_{xx} & \hat{\pmb\Sigma}_{xy} \\ \hat{\pmb\Sigma}_{yx} & \hat{\Sigma}_{yy} \end{array}\right) 
\quad \text{and} \quad 
\hat{\pmb\mu} = \left( \begin{array}{c} \hat{\pmb\mu}_{x} \\ \hat{\mu}_{y} \end{array}\right) 
\end{equation}
\end{linenomath} 
for the joint data  $\{\pmb Z_1, \dots, \pmb Z_n\}$ with $\pmb Z_i = (\pmb X_i^t, Y_i)^t$.

Several authors \citep[see][]{maronna:1986, croux:2003} proposed to achieve robust regression and inference for casewise outliers by robustifying the components in \bref{eq:sol-LS}.  \citet{croux:2003} replaced the empirical covariance matrix and mean by the multivariate S-estimator \citep{davies:1987}. We will refer to this approach as two-step regression (2S-regression). \citet{croux:2003} have shown that under mild assumptions (including symmetry of  $\varepsilon_i$ and  independence of $\varepsilon_i$ and $\pmb X_i$) 2S-regression is Fisher consistent and asymptotically normal even if the S-estimators of multivariate location and scatter themselves are not consistent. Furthermore, 2S-regression is resilient to all kinds of outliers, that is, vertical outliers, bad leverage points, and good leverage points. Note that down-weighting good leverage points could lead to some efficiency loss, but it may also prevent the underestimation of the variance of the estimator, which could be  problematic for inferential purposes \citep[see for example,][]{ruppert:1990}. 

To deal with casewise and cellwise outliers, we propose to use a generalized S-estimator  that uses the consistent filter described in Section \ref{sec:filter}. The estimator is similar to that in \citet{agostinelli:2014}, but with the filter which is consistent for a broader range of distributions. This generality is  needed in the regression setting. Our proposed globally robust regression estimator, called 3S-regression, is given by:
\begin{linenomath} 
\begin{equation}\label{eq:2SGS-reg}
\begin{aligned}
\hat{\pmb\beta}_{3S} &= \hat{\pmb S}_{xx}^{-1} \hat{\pmb S}_{xy} \\
 \hat{\alpha}_{3S} &= \hat{m}_y - \hat{\pmb m}_x^t \hat{\pmb\beta}_{3S}.
\end{aligned}
\end{equation}
\end{linenomath} 
Here, $(\hat{\pmb m}, \hat{\pmb S})$ is a generalized S-estimator computed  as follows:
\begin{enumerate}[leftmargin=*,labelindent=0pt,label=Step \arabic*.]
\item Filter extreme cellwise outliers to prevent cellwise contaminated cases from having large robust Mahalanobis distances in Step 2, and 
\item Down-weight casewise outliers by applying generalized S-estimator (GSE) for multivariate location and scatter \citep{danilov:2012} to the filtered data from Step 1. The GSE is a generalization of the S-estimator for incomplete data that are missing completely at random (MCAR). Since the independent contamination model (ICM) assumes that cells are outlying completely at random, the MCAR assumption is fulfilled if the ICM model holds.
\end{enumerate}

More precisely, consider a  set of  covariates $\{ \pmb X_1, \dots, \pmb X_n\}$. We perform univariate filtering as described in Section \ref{sec:filter} on each variable, $\{X_{1j},\dots,X_{nj}\}$, $j=1,\dots,p$. Let $\{\pmb U_1, \dots,\pmb U_n\}$ be the resulting auxiliary vectors of zeros and ones with zeros indicating the filtered entry in  $\pmb X_i$. More precisely, $\pmb U_i = (U_{i1}, \dots, U_{ip})^t$, where 
\begin{linenomath} 
\[
U_{ij}=  I( \hat\eta^l_{j,n} - \hat s^l_{j,n}\hat t^l_{j,n} \le X_i \le  \hat\eta^u_{j,n} + \hat s^u_{j,n}\hat t^u_{j,n}).
\]
\end{linenomath} 

The goal of the filter is to prevent  propagation of cellwise outliers. If the fraction of cases with at least one flagged cell is very small (below $1\%$, say) then propagation of cellwise outliers is not an issue and the filter can be safely turned off. The procedure that turns  the filter off when the fraction of affected cases is below a given small  threshold, $\xi$,  is considerably simpler to analyze from the asymptotic point of view. Moreover, it retains all the robustness properties derived from the filter.  
Let 
$n_{0} = \# \{1 \le i \le n: \pmb U_i = \pmb 1\} $ be the  number of complete observations after filtering. We set 
\begin{linenomath} 
\begin{equation}\label{eq:filter-switch}
\pmb U^*_i = \pmb 1 I\left(\frac{n-n_0}{n} \le  \xi\right) + \pmb U_i I\left(\frac{n-n_0}{n} >\xi\right), \quad i=1,\dots,n, 
\end{equation}
\end{linenomath} 
with $\xi$ equal to some small threshold.  In this paper we use $\xi=0.01$.

 Finally,  let $\mathbb Z = (\pmb Z_1, \dots, \pmb Z_n)^t$ and  $\mathbb U = ((\pmb U_1^{*}, \dots,\pmb U_n^{*})^t, \pmb 1)$.
 The generalized S-estimator can now be defined as
\begin{linenomath} 
\begin{equation}\label{eq:2SGS}
\begin{aligned}
\hat{\pmb m} &= \hat{\pmb m}_{GS}(\mathbb Z, \mathbb U),\\
\hat{\pmb S} &=\hat{\pmb S}_{GS}(\mathbb Z, \mathbb U ),
\end{aligned}
\end{equation}
\end{linenomath} 
where $\hat{\pmb m}_{GS}$ and $\hat{\pmb S}_{GS}$ are robust multivariate location and scatter generalized S-estimator for incomplete data, $(\mathbb Z, \mathbb U)$, with Tukey's bisquare rho function $\rho_B(t) = \min(1, 1 - (1 - t)^3)$ and 50\% breakdown point \citep[see][for full definition]{danilov:2012}.   
Note that when $\mathbb U = (\pmb 1, \dots, \pmb 1)$ (i.e., when the input data is complete), the generalized S-estimator  reduces to S-estimator \citep{danilov:2012}.

\subsection{Models with continuous and dummy covariates}

For models with continuous and dummy covariates, the direct computation of 3S-regression is likely to fail because the sub-sampling algorithm (needed to compute the generalized S-estimator) is likely to yield collinear subsamples.  In this case, we endow 3S-regression with an iterative algorithm similar to that in \citet{maronna:2000} to deal with continuous and dummy covariates. 

Consider now the following model:
\begin{linenomath} 
\begin{equation}\label{eq:linear-model-dummy}
	Y_i = \alpha  +\pmb X_i^t\pmb\beta_x + \pmb D_i^t\pmb\beta_d + \varepsilon_i
\end{equation}
\end{linenomath} 
for $i=1,\dots,n$ where $\pmb X_i = (X_{i1},\dots,X_{ip_x})^t$ is a $p_x$ dimensional vector of continuous covariates and $\pmb D_i = (D_{i1}, \dots, D_{ip_d})^t$ is a $p_d$ dimensional vector of dummy covariates. 
Set $\mathbb X = (\pmb X_1,\dots,\pmb X_n)^t$, $\mathbb D = (\pmb D_1,\dots,\pmb D_n)^t$, and $\pmb Y = (Y_1, \dots, Y_n)^t$. We assume that the columns in $\mathbb X$ and $\mathbb D$  are linearly independent. 

We  modify the alternating M- and S-regression approach proposed by \citet{maronna:2000}. 
Our algorithm uses 3S-regression to estimate the coefficients of the continuous covariates and regression M-estimators with Huber's rho function $\rho_H(t) = \min(1,  t^2/2)$ \citep{huber:2009} to estimate the coefficients of the dummy covariates. 
More specifically, the algorithm works as follows:
\begin{linenomath} 
\begin{equation}\label{eq:alternating-MS}
\begin{aligned}
(\hat{\alpha}^{(k)}, \hat{\pmb\beta}_x^{(k)}) &= g( \mathbb X, \pmb Y - \mathbb D \hat{\pmb\beta}_d^{(k-1)} ),\\
\hat{\pmb\beta}_d^{(k)} &= M( \mathbb D, \pmb Y - \hat{\alpha}^{(k)}- \hat{\mathbb X}  \hat{\pmb\beta}_x^{(k)}), \quad \text{for} \quad k = 1,\dots, K,
\end{aligned}
\end{equation}
\end{linenomath} 
where $g( \mathbb X, \pmb Y)$ denotes the operation of 3S-regression for a response vector $(\pmb Y, \mathbb X)$ as defined in \bref{eq:2SGS-reg} and $M( \mathbb D, \pmb Y)$ denotes the operation of regression $M$-estimator with no intercept for $(\pmb Y, \mathbb D)$. We let $\hat{\mathbb X}$ be the imputed $\mathbb X$ with the filtered entries imputed by the best linear predictor using $\hat{\pmb m}^{(k)}$ and $\hat{\pmb S}^{(k)}$, the generalized S-estimates at the $k$-th iteration as defined in \bref{eq:2SGS}. We use $\hat{\mathbb X}$ instead of $\mathbb X$ to control the effect of propagation of cellwise outliers.

As  in \citet{maronna:2000}, to calculate the initial estimates, $(\hat{\alpha}^{(0)},  \hat{\pmb\beta}_x^{(0)}, \hat{\pmb\beta}_d^{(0)})$, we first remove  the effect of $\pmb D_i$ from  the continuous covariates and the response variable. Let
\begin{linenomath} 
\[
\overline{\pmb Y} = \pmb Y - \mathbb D \pmb t \quad \text{and} \quad \overline{\mathbb X} = \mathbb X - \mathbb D \mathbb T,
\]
\end{linenomath} 
where $\pmb t = M( \mathbb D, \pmb Y)$ and $\mathbb T$ is a $p_d \times p_x$-matrix with the  $j$-th column as $\pmb T_j = M(\mathbb D, (X_{1j}, \dots, X_{nj})^t)$. Now, the initial estimates are defined by 
\begin{linenomath} 
\[
\begin{aligned}
(\hat{\alpha}^{(0)},  \hat{\pmb\beta}_x^{(0)t}) &= g( \overline{\mathbb X}, \overline{\pmb Y}), \\
\hat{\pmb\beta}_d^{(0)} &=  M( \mathbb D, \pmb Y - \hat{\alpha}^{(0)}- \hat{\mathbb X}  \hat{\pmb\beta}_x^{(0)}).
\end{aligned}
\]
\end{linenomath}

Finally, the procedure in \bref{eq:alternating-MS} is iterated until convergence or until it reaches a maximum of $K=20$ iterations.  
We choose $K = 20$ because our simulation has shown that the procedure usually converges for $K < 20$, provided good initial estimates are used.

\section{Asymptotic properties of three-step regression}\label{sec:asymptotics}

Theorem \ref{thm:3S-S-equiv} (proved in the supplementary material) establishes the equivalence between 3S-regression and 2S-regression \citep{croux:2003} for  the case of continuous covariates.  Let $(\hat \alpha_{3S}, \hat{\pmb\beta}_{3S}^t)$ be the 3S-regression estimate and $(\hat \alpha_{2S}, \hat{\pmb\beta}_{2S}^t)$ be the 2S-regression estimate based on the  sample $\{\pmb Z_1, \dots, \pmb Z_n\}$, where $\pmb Z_i = (\pmb X^t_i, Y_i)$. Let $G(\pmb x)$ and $G_j(x)$ be the distribution functions for $\pmb X_i$  for $X_{ij}$ respectively. 

\begin{theorem}\label{thm:3S-S-equiv}
Suppose that Assumption \ref{assum:tails} holds for each $G_j$, $j=1,\dots,p$. Then, with probability one, for sufficiently large $n$, $\hat \alpha_{3S}= \hat\alpha_{2S}$ and  $\hat{\pmb\beta}_{3S} = \hat{\pmb\beta}_{2S}$. 
\end{theorem}

Since 3S-regression becomes 2S-regression for sufficiently large $n$, 3S-regression inherits  the established asymptotic properties of 2S-regression. Corollary \ref{coro:consistency-normal} states the strong consistency and asymptotic normality of 3S-regression. The corollary requires the following regularity assumptions that are  needed for deriving the consistency and asymptotic normality of 2S-regression \citep[see][]{croux:2003}.

\begin{assumption}\label{assum:error}
Let $F_\varepsilon$ be the distribution of the error term $\varepsilon_i$ in \bref{eq:linear-model}. The distribution $F_\varepsilon$ has a positive, symmetric and unimodal density $f_\varepsilon$.
\end{assumption}

\begin{assumption}\label{assum:carriers}
For all $\pmb v \in \mathbb R^p$ and $\delta \in \mathbb R$,
$P_G( \pmb X_i^t\pmb v = \delta) < 1/2$. 
\end{assumption}

\begin{corollary}\label{coro:consistency-normal}
Suppose that Assumption \ref{assum:tails} holds for each $G_j$, $j=1,\dots,p$, and Assumption \ref{assum:error}--\ref{assum:carriers} hold.  Denote $\hat{\pmb\theta}_{3S} = (\hat\alpha_{3S}, \hat{\pmb\beta}_{3S}^t)^t$ and $\pmb\theta = (\alpha, \pmb\beta^t)^t$. Then, 
\begin{enumerate}[(a)]
\item $\hat\theta_{3S} \to \theta$ a.s..
\item   Let $H$ be the distribution of $(\pmb X^t, Y)$ and let  $(\pmb m_H, \pmb S_H)$ be the S-estimator functional  \citep[see][]{lopuhaa:1989}. We use the same partition outlined in \bref{eq:partition}  for $(\pmb m_H, \pmb S_H)$.
Set $\tilde{\pmb X} = (1, \pmb X^t)^t$.  Then, 
\begin{linenomath} 
\[
\sqrt{n} (\hat{\pmb\theta}_{3S} - \pmb\theta) \to_d  N( \pmb 0, ASV(H) ),
\]
\end{linenomath} 
where 
\begin{linenomath} 
\[
ASV( H) = C(H)^{-1} D(H) C(H)^{-1},
\]
\end{linenomath} 
and where
\begin{linenomath} 
\[
\begin{aligned}
C(H) &= E_{H}\left\{ w( d_{H}(\pmb Z)) \tilde{\pmb X} \tilde{\pmb X}^t\right\} + \frac{2}{\sigma^2_\varepsilon(H)}  E_{H}\left\{ w'( d_{H}(\pmb Z)) (Y - \tilde{\pmb X}^t \pmb\theta)^2 \tilde{\pmb X}\tilde{\pmb X}^t\right\}, \\ 
D(H) &= E_{H}\left\{ w^2( d_H(\pmb Z)) (Y - \tilde{\pmb X}^t \pmb{\theta})^2 \tilde{\pmb X}\tilde{\pmb X}^t \right\}, \\
\sigma_\varepsilon(H) &= \sqrt{S_{H,yy} - \pmb\beta^t \pmb S_{H,xx}\pmb\beta},\\
d_H(\pmb Z) &= (\pmb Z - \pmb m_H)^t \pmb S_H^{-1} (\pmb Z - \pmb m_H), \\
w(t) &= \rho'_B(t).
\end{aligned}
\]
\end{linenomath} 
Here, $\rho_B(t)$ is the Tukey's bisquare rho function.
\end{enumerate} 
\end{corollary}

\begin{remark}
\citet{croux:2003} proved the Fisher consistency of 2S-regression, but the strong consistency also follows from that and  Theorem 3.2 in \citet[][]{lopuhaa:1989}.
\end{remark}

The asymptotic covariance matrix needed for inference can be estimated in the following natural way. 
Let  $(\hat{\pmb m}, \hat{\pmb S})$ be the generalized S-estimate and 
$(\hat\alpha_{3S}, \hat{\pmb\beta}_{3S}^t)$ be the 3S-regression estimate. 
Then, replace $\pmb Z_i = (\pmb X_i^t, Y_i)$ by $\hat{\pmb Z}_i = (\hat{\pmb X}_i^t, Y_i)$ and $\tilde{\pmb X}_i = (1, \pmb X_i^t)^t$ by $\widehat{\tilde{\pmb X}}_i = (1, \hat{\pmb X}_i^t)^t$, where $\hat{\pmb X}_i$ is the best linear prediction of  $\pmb X_i$ (which is possibly incomplete due to filter) using $(\hat{\pmb m}, \hat{\pmb S})$. The identified cellwise outliers in $\pmb X_i$ are filtered and imputed in order to avoid the effect of propagation of outliers on the asymptotic covariance matrix  estimation.
Now, 
\begin{linenomath} 
\[
\widehat{ASV(  H)} = \widehat{C(H)}^{-1} \widehat{D(H)}\widehat{C(H)}^{-1},
\]
\end{linenomath} 
where 
\begin{linenomath} 
\[
\begin{aligned}
\widehat{C(H)} &= \frac{1}{n}\sum_{i=1}^n \left\{ w( d_{n}(\hat{\pmb Z}_i))  + \frac{2}{\hat{\sigma}_{\varepsilon, n}^2}  w'( d_{n}(\hat{\pmb Z}_i))\hat r_i^2 \right\} \widehat{\tilde{\pmb X}}_i \widehat{\tilde{\pmb X}}_i^t, \\
\widehat{D(H)} &= \frac{1}{n}\sum_{i=1}^n  w^2( d_{n}(\hat{\pmb Z}_i)) \hat r_i^2  \widehat{\tilde{\pmb X}}_i \widehat{\tilde{\pmb X}}_i^t, \\
\hat{\sigma}_{\varepsilon, n} &= \sqrt{\hat s_{yy} -\hat{\pmb\beta}_{3S}^t \hat{\pmb S}_{xx}\hat{\pmb\beta}_{3S}},\\
d_n(\hat{\pmb Z}_i) &= (\hat{\pmb Z}_i - \hat{\pmb m})^t \hat{\pmb S}^{-1} (\hat{\pmb Z}_i - \hat{\pmb m}), \\
\hat r_i &= Y_i - \widehat{\tilde{\pmb X}}_i^t\hat{\pmb\theta}_{3S}.
\end{aligned}
\]
\end{linenomath} 
Although the asymptotic covariance matrix formula is valid under clean data, we shall show in Section \ref{sec:simulation} that our proposed inference  remains approximately valid in the presence of a moderate fraction of cellwise and casewise outliers.

In the case of continuous and dummy covariates, \citet{maronna:2000} derived asymptotic results for the alternating regression M- and S-estimates. However, there is no proof of asymptotic results when   regression S-estimators are replaced by 2S-regression. The study of the asymptotic properties of the alternating M- and 2S-regression is  worth of future research.

\section{Simulation}\label{sec:simulation}

We carried out extensive simulation studies in R \citep{R:2015} to investigate the performance of 3S-regression by comparing it with least square (LS)  and two robust alternatives:
\begin{enumerate}[(i)]

\item 2S-regression  as in \citet{croux:2003}. The location and scatter S-estimator with bisquare $\rho$ function and 50\% breakdown point is 
computed by an iterative algorithm that uses an initial MVE estimator. The MVE estimator is computed by sub-sampling with a concentration step. This procedure is implemented in the \verb|R| package \verb|rrcov|, function \verb|CovSest|, option \verb|method="bisquare"| \citep{todorov:2009}; and

\item Shooting S-estimator introduced in \citet{ollerer:2015} with bisquare $\rho$ function and 20\% breakdown point (for each simple regression) as suggested by the authors to attain a good trade-off between robustness and efficiency. The \verb|R| code is available at \url{http://feb.kuleuven.be/Viktoria.Oellerer/software}. 

\end{enumerate}
The generalized S-estimates needed by 3S-regression are computed using the \verb|R| package \verb|GSE|, function \verb|GSE| with default options \citep{leung:2015}. The regression M-estimates needed by the alternating M- and 3S-regression
are computed using the \verb|R| package \verb|MASS|, function \verb|rlm|, option \verb|method="M"| \citep{venables:2002}.

\subsection{Models with continuous covariates}\label{sec:simulation-continuous}

We consider the regression model in \bref{eq:linear-model} with $p = 15$ 
and $n=150,300,500,1000$. 
 The random covariates $\pmb X_i$, $i=1,\dots,n$, are generated from  multivariate normal distribution $N_p(\pmb\mu, \pmb\Sigma)$. We set  $\pmb\mu = \pmb 0$ and $\Sigma_{jj} = 1$ for $j=1,\dots,p$ without loss of generality because GSE in the second step of 3S-regression is location and scale equivariant.
To address the fact that 3S-regression  and the shooting S-estimator are not affine-equivariant,  we consider the random correlation structure for $\pmb\Sigma$ as described in \citet{agostinelli:2014}. We fix the condition number of the random correlation matrix at $100$ to mimic the practical situation for  data sets of similar dimensions. 
Furthermore, to address the fact that the two estimators are not regression equivariant, we randomly generate $\pmb\beta$ as $\pmb\beta= R \,\pmb b$, where $\pmb b$ has a uniform distribution on the unit spherical surface and $R$ is set to $10$.
 We set $\alpha=0$ because GSE is location equivariant.  The response variable $Y_i$ is given by $Y_i = \pmb X_i^t\pmb\beta + \varepsilon_i$, where  
$\varepsilon_i$ are independent (also independent of $\pmb X_i$'s) identically normally distributed with mean 0 and $\sigma = 0.5$. 
Finally, we consider the following scenarios:
\begin{itemize}
\item Clean data: No further changes are done to the data;
\item Cellwise contamination: Randomly replace a fraction $\epsilon$ of the cells in the covariates by outliers $X_{ij}^{cont}=E(X_{ij}) + k \times SD(X_{ij})$ and $\epsilon$ proportion of the responses by outliers $Y_{ij}^{cont}=E(Y_{ij}) + k \times SD(\varepsilon_{i})$, where $k=1,2,\dots,10$;
\item Casewise contamination: Randomly replace a fraction $\epsilon$ of the cases  by leverage outliers $({\pmb X_{i}^{cont}}^t, Y_i^{cont})$, where $\pmb X_{i}^{cont} = c\pmb v$ and $Y_i^{cont} = {\pmb X_i^{cont}}^{t}\pmb\beta + \varepsilon_i^{cont}$ with $\varepsilon_i^{cont} \sim N( k, \sigma^2)$, where $k = 1,2,\dots,15$. Here, $\pmb v$ is the eigenvector corresponding to the smallest eigenvalue of $\pmb\Sigma$ with length such that $(\pmb v - \pmb\mu)^t \pmb\Sigma^{-1}(\pmb v - \pmb\mu) = 1$.  Monte Carlo experiments  show that the placement of outliers in this direction, $\pmb v$, is the least favorable for our estimator. We repeat the  simulation study  in \citet{agostinelli:2014} for dimension 16 and observe that $c = 8$ is the least favorable value for the performance of the scatter estimator.
\end{itemize}
We consider $\epsilon=0.01,0.05$ for cellwise contamination, and $\epsilon=0.10$ for casewise contamination. The number of replicates for each setting is $N = 1000$.

\subsubsection{Coefficient estimation performance}

We examine the effect of cellwise and casewise outliers on the bias of the estimated coefficients. We evaluate the bias using the Monte Carlo mean squared error (MSE):
\begin{linenomath} 
\[
	\overline{MSE}=  \frac{1}{N} \sum_{m=1}^N  \frac{1}{p}\sum_{j=1}^p (\hat\beta_j^{(m)} - \beta_{j}^{(m)} )^2
\]
\end{linenomath} 
where $\hat{\beta}_j^{(m)}$ is the estimate for $\beta^{(m)}_j$ at the $m$-th simulation run.

\begin{table}[t]
\caption{Maximum $\overline{MSE}$ in all the considered scenarios for models with continuous covariates.} 
\label{tab:simulation-maxMSE}
\centering
\small
\begin{tabular}{lrrrrrrrr}
  \hline
 & \multicolumn{2}{c}{Clean} & \multicolumn{2}{c}{1\% Cellwise} & \multicolumn{2}{c}{5\% Cellwise} & \multicolumn{2}{c}{Casewise} \\ 
\cmidrule(lr){2-3} \cmidrule(lr){4-5} \cmidrule(lr){6-7} \cmidrule(lr){8-9}   
$n=$ & 150 & 300 & 150 & 300 & 150 & 300 & 150 & 300 \\  
   \hline 
  3S & 0.012 & 0.005 & 0.039 & 0.020 & 0.902 & 0.797 & 0.223 & 0.143 \\  
  ShootS & 0.034 & 0.017 & 0.134 & 0.080 & 1.129 & 0.912 & 1.570 & 1.460 \\ 
2S & 0.010 & 0.004 & 0.025 & 0.014 & 3.364 & 3.041 & 0.109 & 0.122 \\ 
LS & 0.009 & 0.004 & 2.723 & 2.440 & 4.812 & 4.732 & 8.286 & 8.182 \\ 
   \hline
\end{tabular}
\end{table}

\begin{figure}[t]
\centering
\includegraphics[scale=0.6]{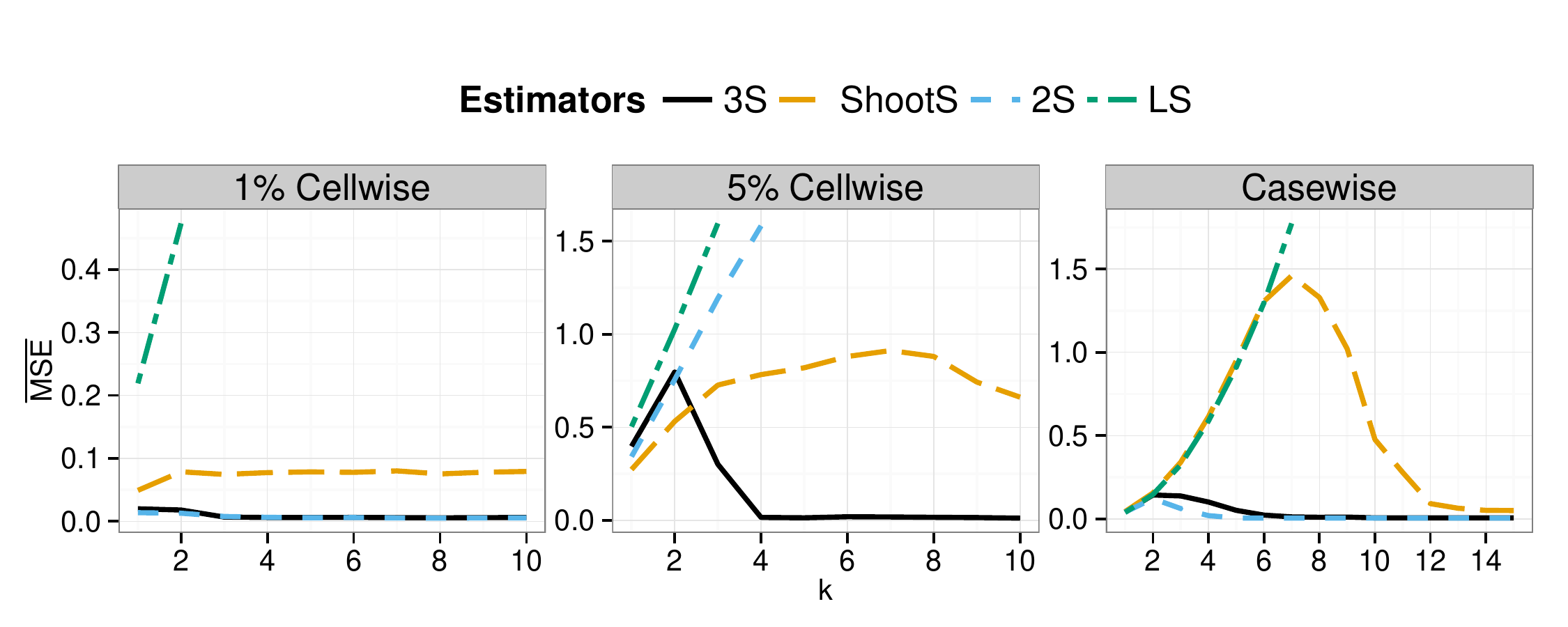}
\caption{$\overline{MSE}$ for various cellwise and casewise contamination values, $k$, for models with continuous covariates. The sample size is $n=300$.}\label{fig:simulation-MSEcurve}
\end{figure}

Table \ref{tab:simulation-maxMSE} shows the $\overline{MSE}$ for clean data  and the  maximum $\overline{MSE}$ for all the cellwise and casewise contamination settings for $n=150,300$. Figure \ref{fig:simulation-MSEcurve} shows the curves of  $\overline{MSE}$ for various cellwise and casewise contamination values for   $n=300$. The results for $n=150$ are similar and the corresponding figure is shown as supplementary material.

In the cellwise contamination setting, 3S-regression is highly robust against moderate and large cellwise outliers ($k \ge 3$), but less robust against inliers ($k  \le  2$). Notice that inliers also affect  the performance of the shooting S-estimator but to a lesser extent. Since the filter does not  flag inliers,  3S-regression and 2S-regression perform similarly in the presence of inliers (see the central panel of Figure \ref{fig:simulation-MSEcurve}).  
The shooting S-estimator is highly robust against large outliers, but less so against moderate cellwise outliers. As expected, 2S-regression  breaks down in the case of $\epsilon = 0.05$, when the propagation of large cellwise outliers is expected to affect more than 50\% of the cases. 

In the casewise contamination setting, 2S-regression has the best performance, as expected. 3S-regression also performs fairly well in this setting. The shooting S-estimator  performs less satisfactorily in this case.

We have also considered other simulation settings and observed similar results (not shown here). In particular, we  considered $p=5$ with $n=50,100$ and $p=25$ with $n=250,500$ under the same set of scenarios (clean data, cellwise contamination, and casewise contamination). Moreover,  we   studied the performance of 3S-regression for larger casewise contamination levels up to 20\%. 3S-regression maintains its competitive performance, outperforming Shooting S and not falling too far behind 2S-regression, which is expected to win in these situations.

\subsubsection{Performance of confidence intervals}

We then assess the performance of confidence intervals for the regression coefficients based on the asymptotic covariance matrix as described in Section \ref{sec:asymptotics}. Intervals that have a coverage close to the nominal
value, while being relatively short, are desirable. 

The $100(1-\tau)\%$ confidence interval (CI) of 3S-regression has the form: 
\begin{linenomath} 
\[
CI(\hat\beta_j) = \left[\hat\beta_j - \Phi^{-1}(1 - \tau/2)\sqrt{\widehat{ASV}(\hat\beta_j)/n}, \,\, \hat\beta_j + \Phi^{-1}(1 - \tau/2)\sqrt{\widehat{ASV}(\hat\beta_j)/n}\right],
\]
\end{linenomath} 
for $j = 0,1, \dots,p$, where $\hat\beta_0 = \hat\alpha$.  We consider $\tau=0.05$ here. We evaluate the performance of CI using the Monte Carlo mean coverage rate (CR): 
\begin{linenomath} 
\[
	\overline{CR} =  \frac{1}{N} \sum_{m=1}^N \frac{1}{p}\sum_{j=1}^p I( \beta_j^{(m)} \in CI(\hat\beta_j^{(m)}) ),
\]
\end{linenomath} 
and the Monte Carlo mean CI lengths:
\begin{linenomath} 
\[
	\overline{CIL} =  \frac{1}{N} \sum_{m=1}^N \frac{1}{p}\sum_{j=1}^p 2\Phi^{-1}(1 - \tau/2)\sqrt{\widehat{ASV}(\hat\beta_j)/n}.
\]
\end{linenomath} 

\begin{figure}[!t]
\centering
\includegraphics[scale=0.6]{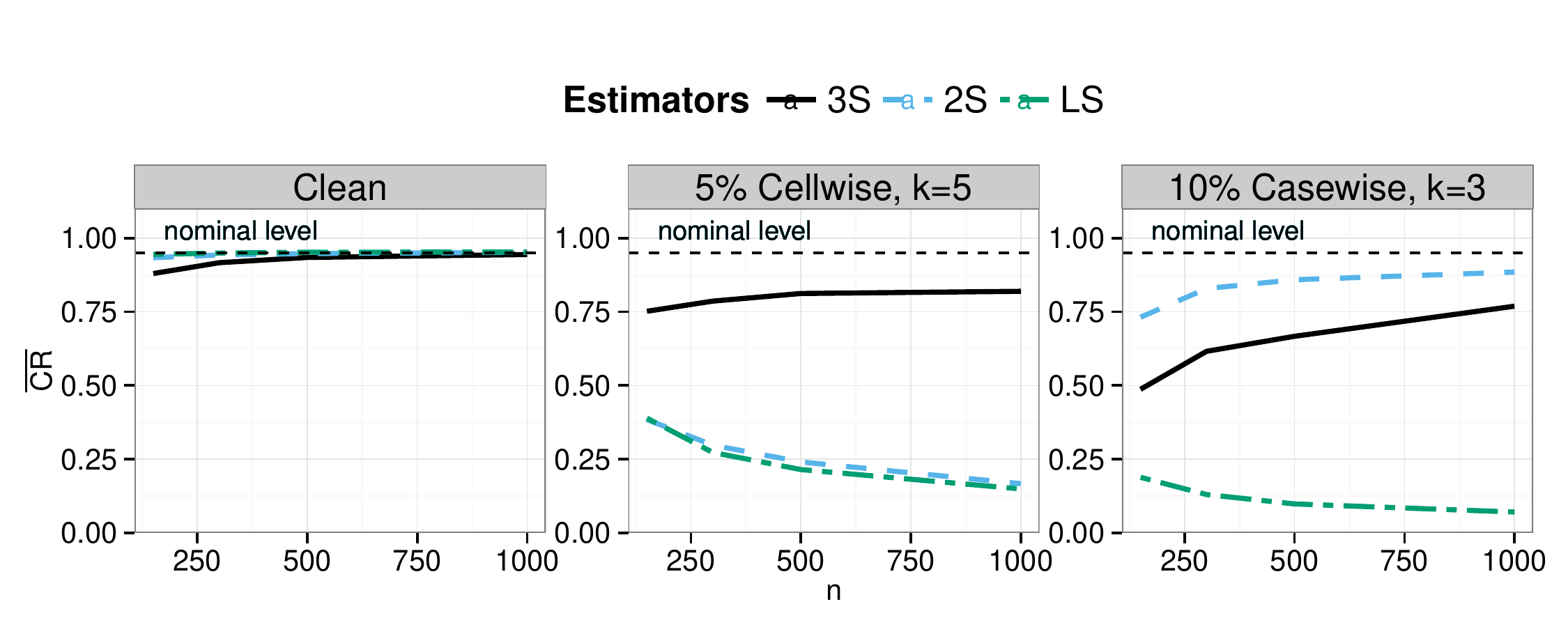}
\caption{$\overline{CR}$ for clean data and for cellwise and casewise contaminated data of various sample size, $n$. }\label{fig:simulation-CI-CoverageCurve}
\end{figure}

Figure \ref{fig:simulation-CI-CoverageCurve}  shows the $\overline{CR}$ in the case of clean data, 5\% cellwise contamination ($k=5$), and 10\% casewise contamination ($k=3$) simulation, with different sample sizes $n=150,300,500,1000$. The nominal value of 95\% is indicated by the horizontal line in the figure.

For clean data, the coverage rates of all the intervals reach the nominal level when the sample size grows, as expected. For data with casewise outliers,  2S-regression yields the best coverage rate, which is closest to the nominal level. However, 3S-regression has an acceptable performance, comparable with that of 2S-regression. For data with cellwise outliers, 3S-regression yields intervals with a coverage rate  relatively closer to the nominal value than LS and 2S-regression.

\begin{table}[!t]
\centering
\caption{Average lengths of confidence intervals for clean data and for cellwise and casewise contamination.}\label{tab:simulation-CI-AverageLength}
\small
\begin{tabular}{rcccccccc}
  \hline
  & \multicolumn{2}{c}{Clean}  & \multicolumn{2}{c}{1\% Cell., $k=5$} & \multicolumn{2}{c}{5\% Cell., $k=5$}  & \multicolumn{2}{c}{10\% Case., $k=3$} \\\cmidrule(lr){2-3} \cmidrule(lr){4-5} \cmidrule(lr){6-7}\cmidrule(lr){8-9}
Size ($n$) & 3S & 2S &  3S & 2S &  3S & 2S &3S & 2S  \\    
  \hline
150 & 0.341 & 0.352 & 0.355 & 0.402 & 0.450 & 1.519 & 0.329 & 0.355 \\ 
   300 & 0.242 & 0.247 & 0.244 & 0.275 & 0.294 & 1.148 & 0.239 & 0.253 \\ 
   500 & 0.187 & 0.189 & 0.190 & 0.212 & 0.222 & 0.912 & 0.189 & 0.197 \\ 
  1000 & 0.133 & 0.133 & 0.134 & 0.150 & 0.155 & 0.662 & 0.137 & 0.140 \\ 
   \hline
\end{tabular}
\end{table}

Furthermore, the length of the intervals obtained from 3S regression is comparable to that LS for clean data and that of 2S-regression for clean data and data with casewise outliers. For data with cellwise outliers, 3S-regression yields intervals with lengths relatively closer to the case of clean data. 
Table \ref{tab:simulation-CI-AverageLength} shows the average lengths of the confidence intervals obtained from 3S- and 2S-regression in the case of clean data, 1\% cellwise contamination ($k=5$), 5\% cellwise contamination ($k=5$), and 10\% casewise contamination ($k=3$) simulation, with different sample sizes $n=150,300,500,1000$. The results of LS are not included here.

In general, 3S-regression yields slightly shorter intervals than 2S-regression in all scenarios because the asymptotic variance is calculated on the data with the filtered cells imputed instead of the complete data. 
On the other hand, 2S-regression tends to yield longer intervals in the cellwise contamination model, even when the propagation of outliers is below the 0.5 breakdown point under THCM, for example, when $\varepsilon = 0.01$. This maybe because 2S-regression loses a significant amount of clean data for estimation when it down-weights cases with outlying components.

\subsection{Models with continuous and dummy covariates}\label{sec:simulation-dummy}

We now conduct a  simulation study to assess the performance of our procedure when  the model includes continuous and dummy covariates. We consider the regression model in \bref{eq:linear-model-dummy} with $p_x = 12$, $p_d = 3$, and $n=150, 300$. 
The random covariates $(\pmb X_i, \pmb D_i)$, $i=1,\dots,n$, are first generated from  multivariate normal distribution $N_p(\pmb 0, \pmb\Sigma)$, where $\pmb{\Sigma}$ is the randomly generated correlation matrix with a fixed condition number of 100. Then,
 we dichotomize $D_{ij}$ at $\Phi^{-1}(\pi_j)$ where $\pi_j = \frac{1}{4},\frac{1}{3},\frac{1}{2}$ for $j=1,2,3$, respectively. Finally, the rest of data are generated in the same way as described in Section~\ref{sec:simulation-continuous}.

In the simulation study, we consider the following scenarios: 
\begin{itemize}
\item Clean data: No further changes are done to the data;
\item Cellwise contamination: Randomly replace a $\epsilon$ fraction of the cells in $\mathbb X$ by outliers $X_{ij}^{cont}=E(X_{ij}) + k \times SD(X_{ij})$ and $\epsilon$ proportion of the responses by outliers $Y_{ij}^{cont}=E(Y_{ij}) + k \times SD(\varepsilon_{i})$, where $k=1,2,\dots, 10$;
\item Casewise contamination:  Let $\pmb\Sigma_x$ be the sub-matrix of $\pmb\Sigma$ with rows and columns corresponding to the continuous covariates.
Randomly replace a $\epsilon$ fraction of the cases in $\mathbb X$ by leverage outliers ${\pmb X_{i}^{cont}} = c\pmb v$, where  $\pmb v$ is the eigenvector corresponding to the smallest eigenvalue of $\pmb\Sigma_x$ with length such that $(\pmb v - \pmb\mu_x)^t \pmb\Sigma^{-1}(\pmb v - \pmb\mu_x) = 1$.  In this case, the number of continuous variables is  13 (instead of 16) and the corresponding least favorable casewise contamination size is found to be $c=7$ (instead of 8) using the same procedure as in Section \ref{sec:simulation-continuous}.
Finally, we replace the corresponding response value by $Y_i^{cont} = {\pmb X_i^{cont}}^{t}\pmb\beta_x + \pmb D_i^t\pmb\beta_d + \varepsilon_i^{cont}$ with $\varepsilon_i^{cont} \sim N( k, \sigma^2)$, where $k = 1,2,\dots,10$. 
\end{itemize}
Again, we consider $\epsilon=0.01,0.05$ for cellwise contamination, and $\epsilon=0.10$ for casewise contamination. The number of replicates for each setting is $N = 1000$.

\begin{table}[!t]
\caption{Maximum $\overline{MSE}$ in all the considered scenarios for models with continuous and dummy covariates.} 
\label{tab:simulation-dummies-MSE}
\centering
\small
\begin{tabular}{lrrrrrrrr}
  \hline
 & \multicolumn{2}{c}{Clean} & \multicolumn{2}{c}{1\% Cellwise} & \multicolumn{2}{c}{5\% Cellwise} & \multicolumn{2}{c}{Casewise} \\ 
\cmidrule(lr){2-3} \cmidrule(lr){4-5} \cmidrule(lr){6-7} \cmidrule(lr){8-9}   
$n=$ & 150 & 300 & 150 & 300 & 150 & 300 & 150 & 300 \\  
   \hline 
  3S & 0.010 & 0.004 & 0.018 & 0.008 & 0.636 & 0.507 & 0.090 & 0.071 \\ 
  ShootS & 0.012 & 0.005 & 0.026 & 0.015 & 0.746 & 0.468 & 0.450 & 0.387 \\ 
  2S & 0.008 & 0.003 & 0.014 & 0.007 & 1.894 & 1.341 & 0.060 & 0.054 \\ 
LS & 0.007 & 0.003 & 2.785 & 2.532 & 5.162 & 4.981 & 1.332 & 1.322 \\ 
   \hline
\end{tabular}
\end{table}

\begin{figure}[ht]
\centering
\includegraphics[scale=0.6]{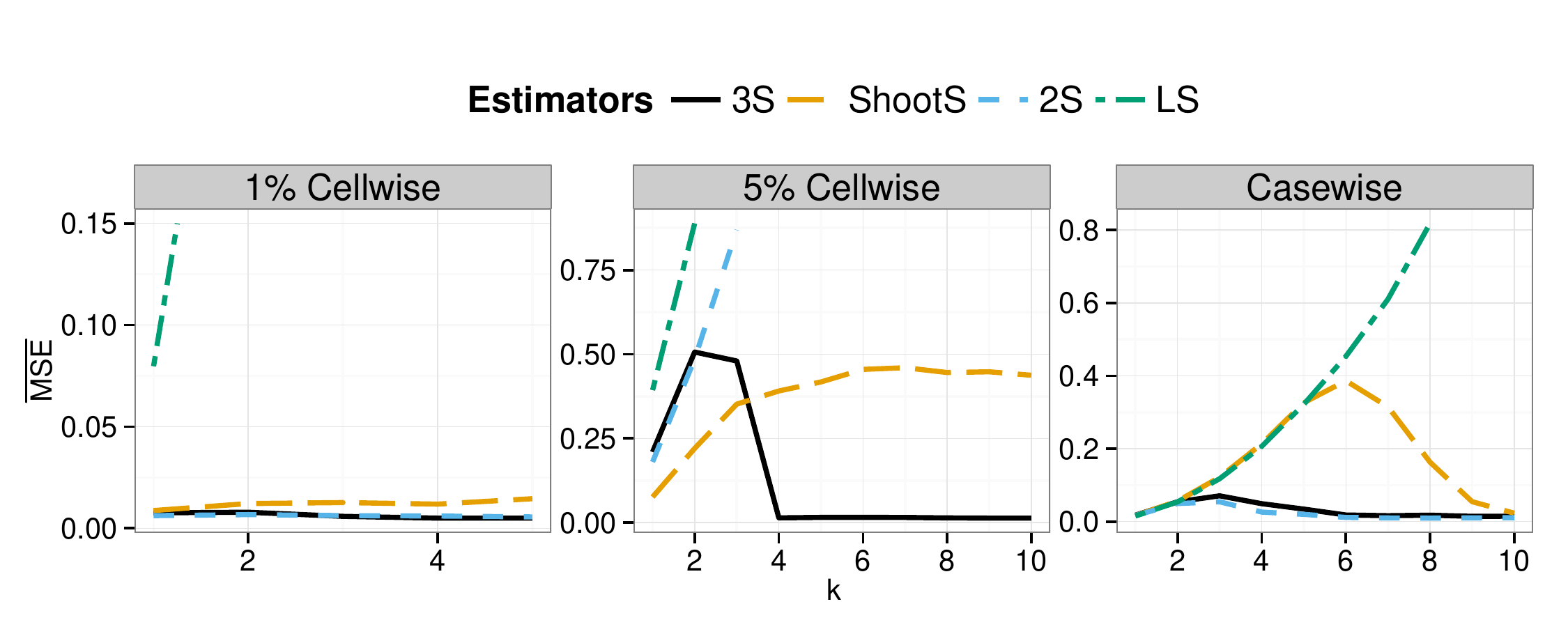}
\caption{$\overline{MSE}$ for various cellwise and casewise contamination values, $k$, for models with continuous and dummy covariates. The sample size is $n=300$.}\label{fig:simulation-dummies-MSEcurve}
\end{figure}

Table \ref{tab:simulation-dummies-MSE} shows the $\overline {MSE}$ for clean data and the  maximum $\overline{MSE}$ for all the cellwise and casewise contamination settings for $n=150,300$.  Figure \ref{fig:simulation-dummies-MSEcurve} shows the curves of  $\overline{MSE}$ for various cellwise and casewise contamination values for   $n=300$. The results for $n=150$ are similar and the corresponding figure is shown as supplementary material. Overall, 3S-regression remains competitive in the case of continuous and dummy covariates.

We also consider the case of non-normal covariates. The covariates are generated from several asymmetric distributions, and the data are contaminated in a similar fashion. The performance of 3S-regression in the case of non-normal covariates is similar to the performance in the case of normal covariates. Results are available as  supplementary material.

\section{Analysis of the Boston housing data}\label{sec:realdata}

We  illustrate the effect of cellwise outlier propagation on  classical robust estimators using the  Boston Housing data. The data, available at the UCI repository \citep{Bache+Lichman:2013},  was collected from  506 census tracts in the Boston Standard Statistical Metropolitan Area in the 1970s on 14 different features.  We consider the nine quantitative variables that were extensively studied \citep[e.g., see in][]{ollerer:2015}. The variables are listed and described in Table 2 in the supplementary material. There is no missing data. The original objective of the study in \citet{harrison:1978} was to analyze the association between the median housing values (medv) in Boston and the residents' willingness to pay for clean air. 

We fit the following model using 3S-regression, the shooting S-estimator, 2S-regression and the LS estimator:
\begin{linenomath} 
\[
\begin{aligned}
log(medv) = \alpha &+ \beta_{1}\, log(crim) + \beta_{2}\, nox^2 + \beta_{3}\, rm^2 + \beta_{x,4}\, age  \\
&+ \beta_{5} \, log(dis) + \beta_{6}\, tax + \beta_{7}\, ptratio + \beta_{8}\, black + \beta_{9}\, log(lstat) + \varepsilon.
\end{aligned}
\]
\end{linenomath} 
The regression coefficient estimates and their P-values are given in Table \ref{tab:boston-estimates}. In particular, we observe that the regression coefficients for the covariates $age$ and $black$ are very different under 3S and 2S-regression. Moreover, $age$  is significant under 2S-regression but highly non-significant under 3S-regression. 
2S-regression is somewhat inefficient  because it throws away a substantial amount of clean data due to the propagation of cellwise outliers. 
It fully down-weights 16.4\% of the  cases in the dataset (cases that receive a zero weight by the multivariate S-estimator). Slightly more than half of these cases (8.7\%) are affected by the propagation of cellwise outliers mainly in the covariates $nox^2$ and $black$  (1.3\% of the cells in the dataset are flagged by the consistent filter). After filtering, these cases have relatively small partial Mahalanobis distances, indicating they are close to the bulk of the data for the remaining variables.

\begin{table}[!t]
\centering
\caption{Estimates and p-values of the regression coefficients for the original Boston Housing data.}\label{tab:boston-estimates}
\footnotesize
\begin{tabular}{lrrrcrrrr}
  \hline
Variable & \multicolumn{2}{c}{3S} & \multicolumn{2}{c}{ShootS} & \multicolumn{2}{c}{2S} & \multicolumn{2}{c}{LS}  \\ 
\cmidrule(lr){2-3}\cmidrule(lr){4-5}\cmidrule(lr){6-7}\cmidrule(lr){8-9}
 & Coeff. & P-Val. & Coeff. & P-Val. & Coeff. & P-Val. & Coeff. & P-Val. \\
  \hline
  log(lstat) & -0.243 & $<$0.001 & -0.266 & - & -0.153 & $<$0.001 & -0.395 & $<$0.001 \\ 
  rm$^2$ & 0.015 & $<$0.001 & 0.013 & - & 0.018 & $<$0.001 & 0.007 & $<$0.001 \\ 
  tax & -0.051 & $<$0.001 & -0.021 & - & -0.046 & $<$0.001 & -0.028 & 0.006 \\ 
  log(dis) & -0.125 & $<$0.001 & -0.157 & - & -0.126 & $<$0.001 & -0.139 & $<$0.001 \\ 
  ptratio & -0.026 & $<$0.001 & -0.027 & - & -0.025 & $<$0.001 & -0.029 & $<$0.001 \\ 
  nox$^2$ & -0.578 & 0.013 & -0.463 & - & -0.445 & 0.023 & -0.451 & $<$0.001 \\ 
  age & -0.023 & 0.645 & -0.040 & - & -0.152 & 0.001 & 0.050 & 0.391 \\ 
  black & -0.726 & 0.398 & 0.787 & - & -0.007 & 0.993 & 0.500 & $<$0.001 \\ 
  log(crim) & -0.006 & 0.513 & 0.004 & - & 0.005 & 0.527 & -0.002 & 0.813 \\ 
\hline
\end{tabular}
\end{table}

We further compare the  four estimators by  computing their squared norm distances, 
$n \times \sum_{j=1}^p (\widehat\beta_{j,A} - \widehat\beta_{j,B})^2 \times MAD( \{ X_{1j}, \dots, X_{nj} \})^2$ \citep[see][]{ollerer:2015}, 
where $MAD$ is the median absolute deviation. Table \ref{tab:boston-pairwise-AND} shows the squared norm distances  for the considered estimators. 
Overall, the three robust estimators are very different from LS. As expected, 3S-regression and shooting S are closer to each other than they are to 2S-regression.  Additional analysis provided as supplementary material indicates  that the observed differences  between the three robust estimators are indeed mostly caused by   the propagation of cellwise outliers in the Boston housing data.

\begin{table}[!t]
\center
\small
\caption{Pairwise squared norm distances between the estimates for the original Boston housing data.}\label{tab:boston-pairwise-AND}
\begin{tabular}{lcccc}
\hline
 & 3S & ShootS & 2S & LS \\ 
\hline
3S & -  & 1.389 & 3.145 & 6.725 \\
ShootS &  & - & 4.312 & 4.661 \\
2S &  &  & - & 16.614 \\
LS &  &  &  & - \\
\hline
\end{tabular}
\end{table}

\section{Concluding remarks}\label{sec:conclusion}

High breakdown point affine equivariant robust estimators are neither efficient nor robust in the independent cellwise contamination model (ICM). By efficiency here  we mean the ability to  use the clean part of the data. In fact, classical robust estimators  are inefficient under  ICM because  they may down-weight an entire row  with a single component being contaminated. Therefore, they may lose some useful information contained in the data. Furthermore, the classical high breakdown point affine equivariant robust estimators may break down under  ICM.  A small fraction  of cellwise outliers  could propagate, affecting a large proportion of cases. For instance, the probability $\overline{\epsilon}$ that at least one component of a case is contaminated is $\overline{\epsilon}=1-(1-\epsilon)^{p}$,
where $\epsilon$ is the proportion of independent cellwise outliers.  This implies that even if $\epsilon$ is small, $\overline{\epsilon}$ could be large for large $p$, and could exceed the 0.5 breakdown point under THCM.  For example, if $\epsilon=0.1$ and
$p=10$, then $\overline{\epsilon}=0.65$; and if $\epsilon=0.05$ and $p=20$, then
$\overline{\epsilon}= 0.64$. 

To overcome these deficiencies of the classical robust estimators, we introduce a three-step regression estimator that can deal with cellwise and casewise outliers. The first step of our estimator is aimed at reducing the impact of outliers propagation posed by  ICM. The second step is aimed at achieving robustness under THCM. As a result, the robust regression estimate from the third step is shown to be efficient (in terms of data usage) and robust under ICM and THCM. We also prove that our estimator is  consistent and asymptotically normal at the central regression model distribution. Finally, we extend our estimator to models with continuous and dummy covariates and provide an algorithm to compute the regression  coefficients. 

The proposed procedures are implemented in the R package \texttt{robreg3S}, which is  freely available on CRAN \citep[the Comprehensive R Archive Network,][]{R:2015}.

\begin{footnotesize}
\section*{Acknowledgement}
\noindent Ruben Zamar's and Andy Leung's research were partially funded by the Natural Science and Engineering Research Council of Canada.

\end{footnotesize}

\bibliographystyle{elsarticle-harv} 
\bibliography{3Sregression}

\begin{thebibliography}{26}
\expandafter\ifx\csname natexlab\endcsname\relax\def\natexlab#1{#1}\fi
\expandafter\ifx\csname url\endcsname\relax
  \def\url#1{\texttt{#1}}\fi
\expandafter\ifx\csname urlprefix\endcsname\relax\def\urlprefix{URL }\fi

\bibitem[{Agostinelli et~al.(2015)Agostinelli, Leung, Yohai, and
  Zamar}]{agostinelli:2014}
Agostinelli, C., Leung, A., Yohai, V.~J., Zamar, R.~H., 2015. Robust estimation
  of multivariate location and scatter in the presence of cellwise and casewise
  contamination. TEST 24~(3), 441--461.

\bibitem[{Alqallaf et~al.(2009)Alqallaf, Van~Aelst, Yohai, and
  Zamar}]{alqallaf:2009}
Alqallaf, F., Van~Aelst, S., Yohai, V.~J., Zamar, R.~H., 2009. Propagation of
  outliers in multivariate data. Ann Statist 37~(1), 311--331.

\bibitem[{Bache and Lichman(2013)}]{Bache+Lichman:2013}
Bache, K., Lichman, M., 2013. {UCI} machine learning repository.
  \url{http://archive.ics.uci.edu/ml}.

\bibitem[{Croux et~al.(2003)Croux, van Aelst, and Dehon}]{croux:2003}
Croux, C., van Aelst, S., Dehon, C., 2003. Bounded influence regression using
  high breakdown scatter matrices. Ann Inst Statist Math 55, 265--285.

\bibitem[{Danilov(2010)}]{danilov:2010}
Danilov, M., 2010. Robust estimation of multivariate scatter under non-affine
  equivarint scenarios. Ph.D. thesis, University of British Columbia.

\bibitem[{Danilov et~al.(2012)Danilov, Yohai, and Zamar}]{danilov:2012}
Danilov, M., Yohai, V.~J., Zamar, R.~H., 2012. Robust estimation of
  multivariate location and scatter in the presence of missing data. J Amer
  Statist Assoc 107, 1178--1186.

\bibitem[{Davies(1987)}]{davies:1987}
Davies, P., 1987. Asymptotic behaviour of {S}-estimators of multivariate
  location parameters and dispersion matrices. Ann Statist 15, 1269--1292.

\bibitem[{Farcomeni(2014{\natexlab{a}})}]{farcomeni:2014a}
Farcomeni, A., 2014{\natexlab{a}}. Robust constrained clustering in presence of
  entry-wise outliers. Technometrics 56, 102--111.

\bibitem[{Farcomeni(2014{\natexlab{b}})}]{farcomeni:2014b}
Farcomeni, A., 2014{\natexlab{b}}. Snipping for robust {K}-means clustering
  under component-wise contamination. Stat Comp 24, 909--917.

\bibitem[{Farcomeni(2015)}]{farcomeni:2015}
Farcomeni, A., 2015. Comments on: Robust estimation of multivariate location
  and scatter in the presence of cellwise and casewise contamination. TEST.

\bibitem[{Fu(1998)}]{fu:1998}
Fu, W., 1998. Penalized regressions: {T}he bridge versus the lasso. J Comput
  Graph Statist 7~(3), 397--416.

\bibitem[{Gervini and Yohai(2002)}]{gervini:2002}
Gervini, D., Yohai, V.~J., 2002. A class of robust and fully efficient
  regression estimators. Ann Statist 30~(2), 583--616.

\bibitem[{Harrison and Rubinfeld(1978)}]{harrison:1978}
Harrison, D., Rubinfeld, D.~L., 1978. Hedonic prices and the demand for clean
  air. J Environ Econ Manage 5, 81--102.

\bibitem[{Huber and Ronchetti(2009)}]{huber:2009}
Huber, P.~J., Ronchetti, E.~M., 2009. Robust Statistics (2nd edition). {John
  Wiley \& Sons}, New Jersey.

\bibitem[{Leung et~al.(2015)Leung, Danilov, Yohai, and Zamar}]{leung:2015}
Leung, A., Danilov, M., Yohai, V., Zamar, R., 2015. GSE: Robust Estimation in
  the Presence of Cellwise and Casewise Contamination and Missing Data. R
  package version 3.2.3.

\bibitem[{Lopuha\"{a}(1989)}]{lopuhaa:1989}
Lopuha\"{a}, H.~P., 1989. On the relation between {S}-estimators and
  {M}-estimators of multivariate location and covariance. Ann Statist 17,
  1662--1683.

\bibitem[{Maronna and Morgenthaler(1986)}]{maronna:1986}
Maronna, R.~A., Morgenthaler, S., 1986. Robust regression through robust
  covariance matrices. Comm Statist Theory Methods 15, 1347--1365.

\bibitem[{Maronna and Yohai(2000)}]{maronna:2000}
Maronna, R.~A., Yohai, V.~J., 2000. Robust regression with both continuous and
  categorical predictors. J Statist Plann Inference 89, 197--214.

\bibitem[{{\"O}llerer et~al.(2015){\"O}llerer, Alfons, and
  Croux}]{ollerer:2015}
{\"O}llerer, V., Alfons, A., Croux, C., 2015. The shooting {S}-estimator for
  robust regression. Comput Statist.
\newline\urlprefix\url{http://dx.doi.org/10.1007/s00180-015-0593-7}

\bibitem[{{R Core Team}(2015)}]{R:2015}
{R Core Team}, 2015. R: A Language and Environment for Statistical Computing. R
  Foundation for Statistical Computing, Vienna, Austria.
\newline\urlprefix\url{http://www.R-project.org/}

\bibitem[{Rousseeuw(1984)}]{rousseeuw:1984b}
Rousseeuw, P., 1984. Least median of squares regression. J Amer Statist Assoc
  79, 871--880.

\bibitem[{Rousseeuw and Yohai(1984)}]{rousseeuw:1984a}
Rousseeuw, P.~J., Yohai, V.~J., 1984. Robust regression by means of
  {S}-estimators. In: Franke, J., H{\"a}rdle, W., Martin, D. (Eds.), Robust and
  Nonlinear Time Series. Vol.~26 of Lecture Notes in Statistics. Springer, New
  York, US, pp. 256--272.

\bibitem[{Ruppert and Simpson(1990)}]{ruppert:1990}
Ruppert, D., Simpson, D., 1990. Unmasking multivariate outliers and leverage
  points: Comment. J Amer Statist Assoc 85, 644--646.

\bibitem[{Todorov and Filzmoser(2009)}]{todorov:2009}
Todorov, V., Filzmoser, P., 2009. An object-oriented framework for robust
  multivariate analysis. Journal of Statistical Software 32~(3), 1--47.
\newline\urlprefix\url{http://www.jstatsoft.org/v32/i03/}

\bibitem[{Venables and Ripley(2002)}]{venables:2002}
Venables, W.~N., Ripley, B.~D., 2002. Modern Applied Statistics with S, 4th
  Edition. Springer, New York, iSBN 0-387-95457-0.
\newline\urlprefix\url{http://www.stats.ox.ac.uk/pub/MASS4}

\bibitem[{Yohai(1985)}]{yohai:1985}
Yohai, V.~J., 1985. High breakdown point and high efficiency robust estimates
  for regression. Tech. Rep.~66, Department of Statistics, University of
  Washington, available at
  \url{http://www.stat.washington.edu/research/reports/1985/tr066.pdf}.

\end{thebibliography}


\begin{thebibliography}{2}
\expandafter\ifx\csname natexlab\endcsname\relax\def\natexlab#1{#1}\fi
\expandafter\ifx\csname url\endcsname\relax
  \def\url#1{\texttt{#1}}\fi
\expandafter\ifx\csname urlprefix\endcsname\relax\def\urlprefix{URL }\fi

\bibitem[{Danilov et~al.(2012)Danilov, Yohai, and Zamar}]{danilov:2012}
Danilov, M., Yohai, V.~J., Zamar, R.~H., 2012. Robust estimation of
  multivariate location and scatter in the presence of missing data. J Amer
  Statist Assoc 107, 1178--1186.

\bibitem[{Yohai(1985)}]{yohai:1985}
Yohai, V.~J., 1985. High breakdown point and high efficiency robust estimates
  for regression. Tech. Rep.~66, Department of Statistics, University of
  Washington, available at
  \url{http://www.stat.washington.edu/research/reports/1985/tr066.pdf}.

\end{thebibliography}

\end{document}


\title{Supplementary Material for: Robust regression estimation and inference in the presence of cellwise and casewise contamination}
\author{{\bf Andy Leung}}
\author{{\bf Hongyang Zhang}}
\author{{\bf Ruben H. Zamar}}
\affil{{\small Department of Statistics, University of British Columbia, 3182-2207 Main Mall, Vancouver, British Columbia V6T 1Z4, Canada}}
\renewcommand\Authands{, }

\maketitle

\begin{abstract}
This supplementary material contains all the proofs, additional simulation results, and related supplementary material referenced in the article ``Robust regression estimation and inference in the presence of cellwise and casewise contamination". 
\end{abstract}

\section{Proofs of Lemmas and Theorems}

\subsection{Proof of Theorem 2.1}

We need to following lemma in the proof.

\begin{lemma}\label{lem:consistent-scale}
Let $X, X_1, \dots, X_n$ be independent with a continuous distribution function $G(x)$. Given $0 < \alpha < 1$, let $\eta = G^{-1}(1 - \alpha)$ and $s = \text{med}(X - \eta | X > \eta)$. Now, consider the following estimator: $\hat\eta_n = \hat G_n^{-1}(1 - \alpha)$ and $\hat s_n = \text{med}(\{ X_i - \hat\eta_n | X_i > \hat\eta_n\})$. Then, $\hat s_n \to s$ a.s.
\end{lemma}

\begin{proof}[\bf{Proof}]
Without loss of generality, assume that $X_1 < X_2 < \dots < X_n$. So, $\hat\eta_n = \hat G_n^{-1}(1-\alpha) = X_{\lceil n(1-\alpha)\rceil}$, and 
\begin{linenomath} 
\[
\#\{ X_i | X_i > \hat\eta_n\} = n - \lceil n(1-\alpha)\rceil = n - (n + \lceil -n\alpha \rceil ) = \lfloor n\alpha \rfloor. 
\]
\end{linenomath} 
Then, $X_k = \text{med}(\{ X_i | X_i > \hat\eta_n\})$ where
\begin{linenomath} 
\[
\begin{aligned}
k &= \lceil n(1-\alpha)\rceil + \left\lceil\frac{\lfloor n\alpha \rfloor}{2}\right\rceil \\
&= n - \lfloor n\alpha \rfloor + \left\lceil\frac{\lfloor n\alpha \rfloor}{2}\right\rceil\\
&= n - \left\lfloor\frac{\lfloor n\alpha \rfloor}{2}\right\rfloor = n - \left\lfloor\frac{ n\alpha }{2}\right\rfloor.
\end{aligned}
\]
\end{linenomath} 
In other words, $\text{med}(\{ X_i | X_i > \hat\eta_n\})= X_{\lceil n(1-\alpha/2)\rceil} = \hat G_n^{-1}( 1 - \alpha/2)$, and 
\begin{linenomath} 
\[
\hat s_n = \hat G_n^{-1}( 1 - \alpha/2) - \hat\eta_n.
\]
\end{linenomath} 
Therefore, $\hat s_n \to s$ a.s., where $s = G^{-1}(1 - \alpha/2) - \eta$. 

\end{proof}

\begin{proof}[\bf{Proof of Theorem 2.1}]

Without loss of generality, we consider only the upper tail. Also, to simplify the notation, we drop out the $G$ in the probability and the $u$ that was used to distinguish between the notations for upper tail and lower tail.

Define $F(t)$ and $\hat F_{n}(t)$ by 
\begin{linenomath} 
\[
\begin{aligned}
F(t) = \frac{ P( 0 < ( X - \eta)/s \le t) }{ P( X > \eta) } \quad \text{and}\quad 
\hat F_{n}(t) = \frac{\frac{1}{n}\sum_{i=1}^n I(0 <  (X_i - \hat\eta_n)/\hat s_n \le t) }{ \frac{1}{n} \sum_{i=1}^n I(X_i > \hat\eta_n) }.
\end{aligned}
\]
\end{linenomath} 
Let $F_0(t) = 1 - e^{-t}$. 
It is sufficient to prove that for every $\epsilon > 0$ there exists $N$ such that for all $n \ge N$,
\begin{linenomath} 
\[
\sup_{t \ge t_0} \left\{ F_0(t) - \hat F_n(t) \right\}^+ < \epsilon.
\]
\end{linenomath} 
Note that 
\begin{linenomath} 
\[
\begin{aligned}
|F_0(t) - \hat F_n(t) | &\le  \left|F_0(t) -\frac{ P( 0 < (  X - \eta)/s \le t) }{ P( X \ge \eta) } \right| \\
&\quad + \left| \frac{ P( 0 < (  X - \eta)/s \le t) }{ P( X > \eta) }  - \frac{ P( 0 < (  X -\hat\eta_n)/\hat s_n \le t) }{ P( X > \eta) } \right|\\
&\quad + \left|  \frac{ P( 0 < (  X -\hat\eta_n)/\hat s_n \le t) }{ P( X > \eta) }  - \frac{ P( 0 < (  X -\hat\eta_n)/\hat s_n \le t) }{ P( X > \hat\eta_n) }  \right|\\
&\quad +  \left|  \frac{ P( 0 < (  X -\hat\eta_n)/\hat s_n \le t) }{ P(  X > \hat\eta_n) } -   \frac{ \frac{1}{n}\sum_{i=1}^n I(0 < ( X_i - \hat\eta_n)/\hat s_n \le t) }{ P( X > \hat\eta_n) } \right|\\
&\quad + \left|   \frac{ \frac{1}{n}\sum_{i=1}^n I(0 < ( X_i - \hat\eta_n)/\hat s_n \le t) }{ P( X > \hat\eta_n) } - \frac{\frac{1}{n}\sum_{i=1}^n I(0 <  ( X_i - \hat\eta_n)/\hat s_n \le t) }{ \frac{1}{n} \sum_{i=1}^n I( X_i > \hat\eta_n) } \right| \\
& \quad = A + B + C + D + E.
\end{aligned}
\]
\end{linenomath} 

By Assumption 2.1, $A = 0$.

Note that 
\begin{linenomath} 
\[
\begin{aligned}
B  &= \frac{1}{\alpha} |P( 0 < (  X - \eta)/s \le t) -  P( 0 < (  X -\hat\eta_n)/\hat s_n \le t) |\\
&= \frac{1}{\alpha} |[G(st + \eta) - G(\hat s_n t + \hat\eta_n)] - [G(\eta) - G(\hat\eta_n)] |.
\end{aligned}
\]
\end{linenomath} 
Next, we show that $\sup_t |G(st + \eta) - G(\hat s_n t + \hat\eta_n)| < \varepsilon\alpha/4$ and $|G(\eta) - G(\hat\eta_n)| < \varepsilon\alpha/4$. 

Given a small $\delta_0 > 0$ such that $s - \delta_0 > c$ and $\eta - \delta_0 > c$ for $c > 0$. Choose a large $K > 0$ such that for $K_{\delta_0} = (s - \delta_0)K + \eta - \delta_0$, $G(K_{\delta_0}) > 1 - \frac{\varepsilon\alpha}{4}$. First, consider $t > K$. 
Since $\delta_0 > 0$, we have 
$st + \eta > (s - \delta_0)K+ (\eta - \delta_0) = K_{\delta_0}$, and therefore, $G(st + \eta)> G(K_{\delta_0})$.
Also, by Lemma \ref{lem:consistent-scale}, $\hat s_n \to s$ a.s. and $\hat\eta_{n} \to \eta$ a.s.. So,  there exists $N_0$ such that $|\hat s_n - s|< \delta_0$ and $|\hat\eta_{n} - \eta| < \delta_0$ for all $n \ge N_0$. So, we have $\hat s_n > s - \delta_0$ and $\hat\eta_n > \eta - \delta_0$, which implies
$\hat s_n t + \hat\eta_{n} > (s - \delta_0) K + (\eta - \delta_0) = K_{\delta_0}$ and $G(\hat s_n t + \hat\eta_n) > G(K_{\delta_0})$. Therefore, 
\begin{linenomath} 
\[
\sup_{t > K} \{G(st + \eta) - G(\hat s_n t + \hat\eta_{n}) \} \le \frac{\varepsilon \alpha}{4}.
\]
\end{linenomath} 

Now, consider $t \le K$. We have $|(st  + \eta) - (\hat s_n t + \hat\eta ) | \le t |s - \hat s_n| + |\eta - \hat\eta_{n}| < K \delta_0 + \delta_0 < \delta_1$. 
Now by the uniform continuity of $G$, given $\varepsilon > 0$, there exists $N_1$ such that for $n \ge N_1$, $|(st + \eta) - (\hat s_n t + \hat\eta_{n}) | < \delta$, and therefore, $|G(st + \eta) - G(\hat s_n t + \hat\eta_{n})| < \frac{\varepsilon \alpha}{4}$. Similarly,  there exists $N_2$ such that for $n \ge N_2$, $|\eta - \hat\eta_{n} | < \delta$, and therefore,  $| G(\eta) - G(\hat\eta_{n}) | < \frac{\varepsilon \alpha}{4}$. 

So, with probability one, take  $N = \max\{N_0, N_1, N_2\}$ such that for $n \ge N$, it  implies that $|(st + \eta) - (\hat s_nt  + \hat\eta_{n}) | < \delta$ and $|\eta - \hat\eta_{n}| < \delta$. Then, we have
\begin{linenomath} 
\[
\begin{aligned}
B &\le \frac{1}{\alpha} \left[ \sup_t  \left|G(st + \eta) - G(\hat s_n t + \hat\eta_{n})\right| + \left| G(\eta) - G(\hat\eta_{n}) \right| \right] \\
&\le \frac{1}{\alpha}( \frac{\varepsilon\alpha}{4} + \frac{\varepsilon\alpha}{4} ) = \frac{\varepsilon}{2}.
\end{aligned}
\]
\end{linenomath} 

Next, we have 
\begin{linenomath} 
\[
\begin{aligned}
C &\le \frac{|G(\hat s_n t + \hat\eta_{n}) - G(\hat\eta_{n})|}{(1 - G(\eta))(1 - G(\hat\eta_{n}))}  |G(\eta) - G(\hat\eta_{n}) | \\
&\le \frac{1}{1 - G(\eta)}  |G(\eta) - G(\hat\eta_{n}) | \le \frac{1}{\alpha} \frac{\varepsilon\alpha}{4} =  \frac{\varepsilon}{4}.
\end{aligned}
\]
\end{linenomath} 

By the Gilvenko--Cantelli Theorem, with probability one, we can show that there exists $N_3$ such that for $n \ge N_3$, 
$\sup_t |P(0 < (X - \hat\eta_{n})/\hat s_n \le t) - \frac{1}{n}\sum_{i=1}^n I(0 < (X - \hat\eta_{n})/\hat s_n \le t) |  < \frac{\varepsilon\alpha}{16}.$
Note that for large enough $n$, we have 
 $P(X > \hat\eta_{n}) > \frac{\alpha}{2}$.
So,
\begin{linenomath} 
\[
D = \left|  \frac{ P( 0 < ( X -\hat\eta_{n})/\hat s_n \le t) }{ P( X > \hat\eta_{n}) } -   \frac{ \frac{1}{n}\sum_{i=1}^n I(0 <  (X_i - \hat\eta_{n})/\hat s_n \le t) }{ P( X > \hat\eta_{n}) } \right| \le \frac{2}{\alpha} \frac{\varepsilon\alpha}{16} = \frac{\varepsilon}{8}.
\]
\end{linenomath} 

Next, by the Gilvanko--Cantelli Theorem again, there exists $N_4$ such that for $n \ge N_4$, $|P(X > \hat\eta_{n}) - \frac{1}{n}\sum_{i=1}^n I( X_i > \hat\eta_{n})| < \sup_t |P(X > t) - \frac{1}{n}\sum_{i=1}^n I( X_i > t)| < \frac{\varepsilon\alpha}{16}$. Then, we have
\begin{linenomath} 
\[
\begin{aligned}
E &\le (\frac{1}{n}\sum_{i=1}^n I(0 < (X_i - \hat\eta_{n})/\hat s_n \le t)) \frac{\left|P(X > \hat\eta_{n}) - \frac{1}{n}\sum_{i=1}^n I( X_i > \hat\eta_{n})\right|}{P(X > \hat\eta_{n}) \, (\frac{1}{n}\sum_{i=1}^n I( X_i > \hat\eta_{n}) )}\\
&\le (\frac{1}{n}\sum_{i=1}^n I( X_i > \hat\eta_{n}))\frac{\left|P(X > \hat\eta_{n}) - \frac{1}{n}\sum_{i=1}^n I( X_i > \hat\eta_{n})\right|}{P(X > \hat\eta_{n}) \, (\frac{1}{n}\sum_{i=1}^n I( X_i > \hat\eta_{n}))} \\
& \le \frac{2}{\alpha} \frac{\varepsilon\alpha}{16} = \frac{\varepsilon}{8}.
\end{aligned}
\]
\end{linenomath} 
Finally, take $N = \max\{ N_0, N_1, N_2, N_3, N_4\}$, we have
\begin{linenomath} 
\[
\sup_t \{F(t) - \hat F_n(t) \} \le A + B + C + D + E\le \frac{\varepsilon}{2} + \frac{\varepsilon}{4} + \frac{\varepsilon}{8} + \frac{\varepsilon}{8} = \varepsilon.
\]
\end{linenomath} 
\end{proof}

\subsection{Proof of Theorem 4.1}

 Let $(\pmb U_1, \dots, \pmb U_n)^t$ be the matrix of zeros and ones with zero corresponding to a filtered component in $(\pmb X_1, \dots, \pmb X_n)^t$, and $n_0$ be the number of complete observations after the filter step. Now, let $C_j = \{ i, 1 \le i \le n : U_{ij} = 0\}$ and $C = \cup_{j=1}^p C_j$. So, $C_j$ is the set of indices of filtered values for variable $j$, and $C$ is the set of indices of incomplete observations. By Boole's inequality, 
\begin{linenomath} 
\[
n - n_0 = \# C \le \sum_{j=1}^p \# C_j.
\]
\end{linenomath} 
Let $\xi > 0$ be as described in Section 3. Now, for each variable $\{X_{1j}, \dots, X_{nj}\}$, $j =1, \dots, p$, apply Theorem 2.1 to obtain $N_j$ such that, with probability one, 
\begin{linenomath} 
\[
\# C_j \le n\xi/p, \quad \text{for} \quad n \ge N_j.
\] 
\end{linenomath} 
Set $N = \max\{ N_1, \dots, N_p\}$. Hence, with probability one, 
\begin{linenomath} 
\[
n - n_0 \le  \sum_{j=1}^p \# C_j \le\sum_{j=1}^p n\xi/p = n\xi,
\]
\end{linenomath} 
for $n \ge N$, or equivalently,
\begin{linenomath} 
\[
\begin{aligned}
 \frac{n_0}{n} \ge 1 - \xi.
\end{aligned}
\]
\end{linenomath} 
Therefore, $\pmb U^*_i = (1, \dots, 1)^t$ according to (6), and 
$\mathbb U = \mathbb I$, where $\mathbb I$ has every entry  equal to 1. In other words, for $n \ge N$, the GSE in Section 3 becomes 
\begin{linenomath} 
\[
\begin{aligned}
\hat{\pmb m} &= \hat{\pmb m}_{GS}(\mathbb Z, \mathbb I)\\
\hat{\pmb S} &=\hat{\pmb S}_{GS}(\mathbb Z, \mathbb I ).
\end{aligned}
\]
\end{linenomath} 
Since GSE on complete data reduces to the regular S-estimator \citep{danilov:2012}, this implies that 3S-regression reduces to S-regression for $n \ge N$.

\newpage

\section{Additional figures from the simulation study in Section 5}

\begin{figure}[h!]
\centering
\includegraphics[scale=0.6]{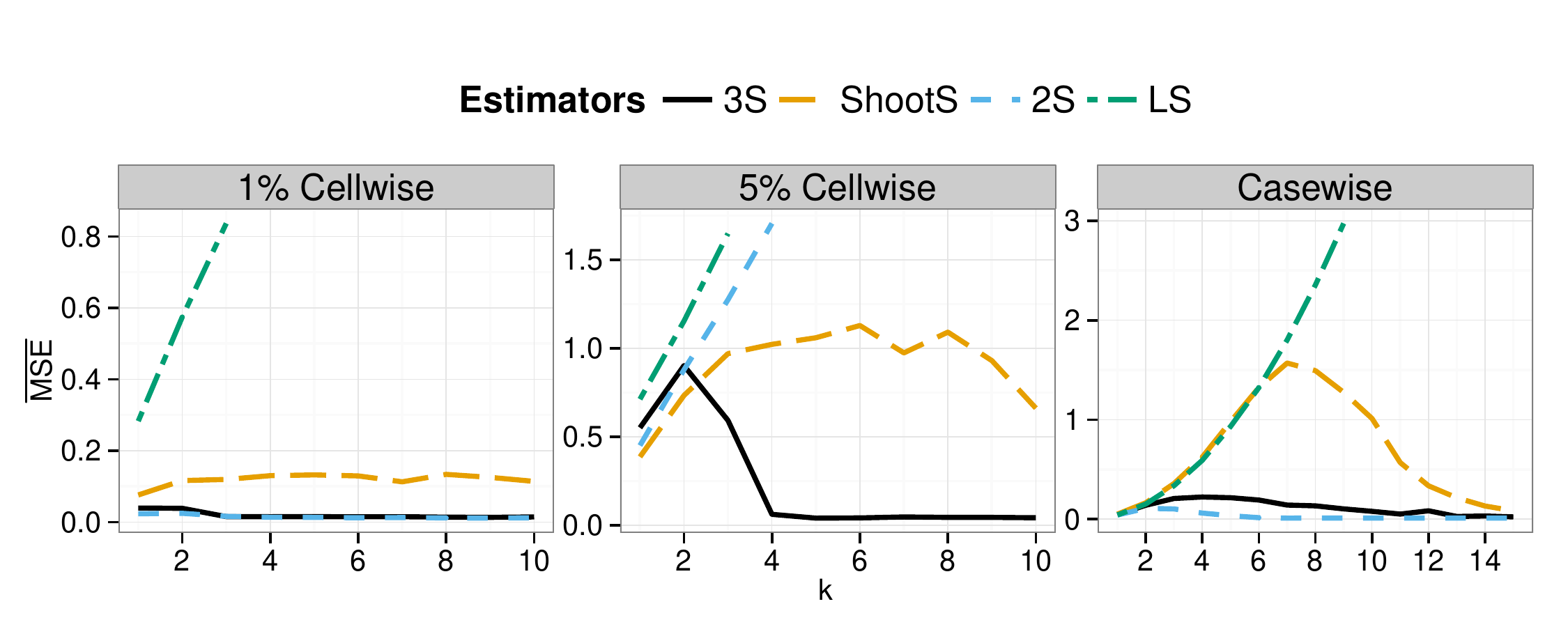}
\caption{$\overline{MSE}$ for various cellwise and casewise contamination values, $k$, for models with $p=15$ continuous covariates. The sample size is $n=150$. For details see Section 5.1 in the paper. }
\end{figure}

\begin{figure}[h!]
\centering
\includegraphics[scale=0.6]{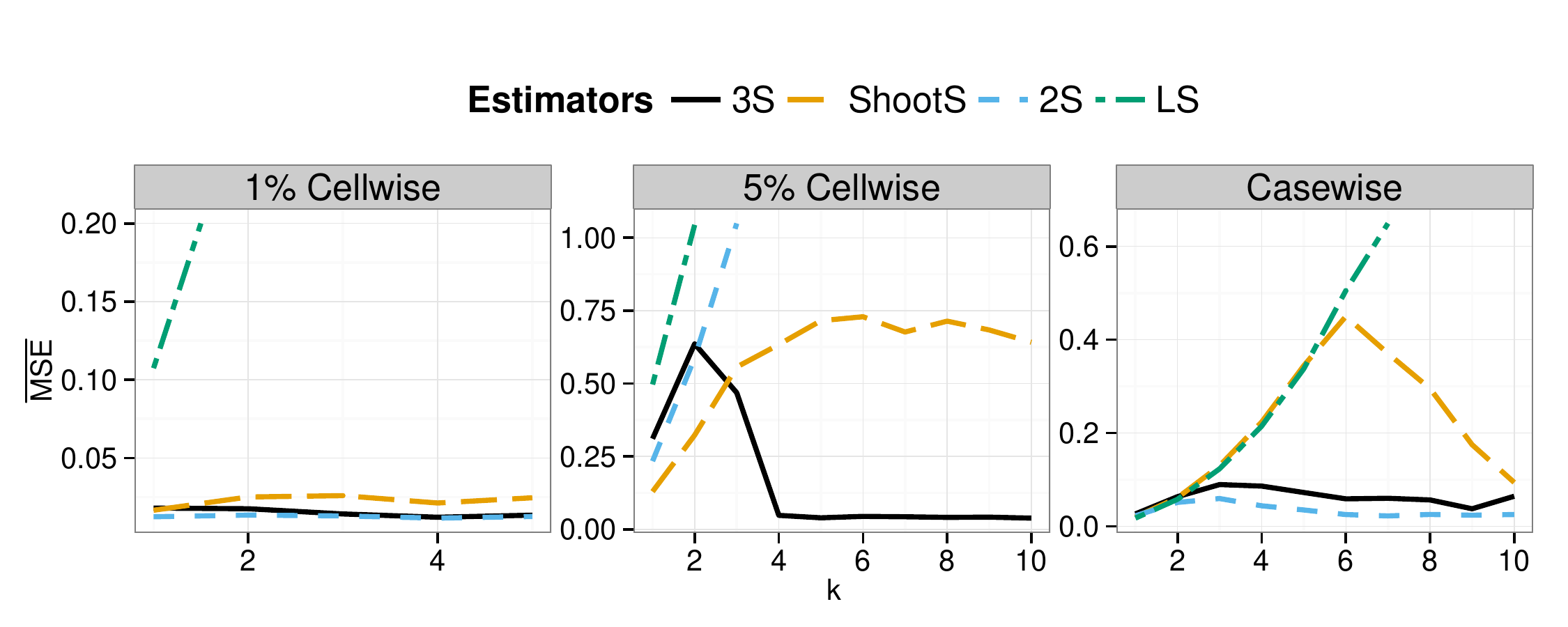}
\caption{$\overline{MSE}$ for various cellwise and casewise contamination values, $k$, for models with $p_x =12$ continuous and $p_d = 3$ dummy covariates. The sample size is $n=150$. For details see Section 5.2 in the paper.}
\end{figure}

\section{Investigation on the performance on non-normal covariates}

Here, we conduct a modest simulation study to compare the performance of 3S-regression, the shooting S-estimator, 2S-regression and the LS estimator for data with non-normal covariates. 

We consider the same regression model with $p = 15$ and $n=300$ as in Section 5, but the covariates are generated from a non-normal distribution as follows. The random covariates $\pmb X_i$, $i=1,\dots,n$, are first generated from  multivariate normal distribution $N_p(\pmb 0, \pmb\Sigma)$, where $\pmb{\Sigma}$ is the randomly generated correlation matrix with a fix condition number of 100. Then, we transform the variables by doing the following:
\begin{linenomath} 
\[
(X_{i1}, X_{i2}, \dots, X_{ip}) \gets (G_1^{-1}(\Phi(X_{i1})), G_2^{-1}(\Phi(X_{i2})), \dots, G_p^{-1}(\Phi(X_{ip})) ),
\]
\end{linenomath} 
where $\Phi(x)$ is the standard normal. We set $G_j$ as $N(0,1)$ for $j=1,2,3$, $\chi^2(20)$ for $j=4,5,6$, $F(90, 10)$ for $j=7,8,9$, $\chi^2(1)$ for $j=10,11,12$, and Pareto$(1,3)$ for $j=13,14,15$.

In the simulation study, we consider the following scenarios:
\begin{itemize}
\item Clean data: No further changes are done to the data;
\item Cellwise contamination: Randomly replace $\epsilon=0.05$ fraction of the cells in the covariates by outliers $X_{ij}^{cont}=k \times G_j^{-1}(0.999)$ and $\epsilon$ proportion of the responses by outliers $Y_{ij}^{cont}=E(Y_{ij}) + k \times SD(\varepsilon_{i})$. We present the results for $k=1,5,10$, but for larger values of $k$ we obtain similar results. 
\end{itemize}
The number of replicates for each setting is $N = 1000$.

The performance of the estimator in terms of $\overline {MSE}$ are summarized in Table \ref{tab:simulation-nonnormal-MSE}. The performance of 3S-regression is comparable  to that of LS and 2S-regression for clean data and outperforms the shooting S, LS and 2S-regression for cellwise-contaminated data, even under some deviations from the assumptions on the  tail distributions of the covariates.   

\begin{table}[thb]
\centering
\small
\caption{$\overline{MSE}$ for clean data and cell-wise contaminated data.}\label{tab:simulation-nonnormal-MSE}
\begin{tabular}{lcccc}
  \hline
Estimators & Clean & \multicolumn{3}{c}{Cellwise}  \\ \cmidrule(lr){3-5}
 & & $k=2$ & $k=5$ & $k=10$ \\
\hline
  3S & 0.007 & 0.014 & 0.013 & 0.015 \\ 
  ShootS & 0.254 & 0.839 & 1.048 & 0.882 \\ 
2S & 0.003 & 4.102 & 3.851 & 4.057 \\ 
LS & 0.001 & 4.311 & 6.438 & 6.588 \\ 
   \hline
\end{tabular}
\end{table}

\section{Further analysis of the Boston housing data}

\begin{table}[!h]
\caption{Description of the variables in the Boston Housing data}\label{tab:boston-variables}\label{tab:boston-variables}
\centering
\begin{footnotesize}
\begin{tabular}{ll}
\hline
Variables & Description \\
\hline
medv (response) & corrected median value of owner-occupied homes in USD 1000's \\
crim   & per capita crime rate by town \\
nox	& nitric oxides concentration (parts per 10 million)\\
rm	& average number of rooms per dwelling\\
age	& proportion of owner-occupied units built prior to 1940\\
dis	& weighted distances to five Boston employment centers\\
tax	& full-value property-tax rate per USD 1,000,000\\
ptratio	& pupil-teacher ratio by town\\
black	& $(B - 0.63)^2$ where B is the proportion of blacks by town\\
lstat	& percentage of lower status of the population\\
\hline
\end{tabular}
\end{footnotesize}
\end{table}

We now further illustrate how  the propagation of cellwise outliers in the Boston housing data leads to the observed differences  among the three robust estimators.

Recall that half of the cases fully downweighted by 2S-regression have entries flagged as cellwise outliers. We  replace these flagged cells  by their best linear predictions (using the 3S-regression estimate) and then, refit the model with the four considered estimators. The resulting coefficient estimates and their P-values   are given  in Table \ref{boston-estimates-imputedver}. Notice that the covariate $age$ is no longer significant under 2S-regression. Moreover, Table \ref{tab:boston-pairwise-AND-experiment}  shows the norm distances between all the estimates calculated from such imputed data. Now,  2S-regression is  considerably closer  to the  cellwise robust estimators, and it no longer fully down-weights the  cases formerly affected by cellwise outliers (the median weight of these cases is now 0.64, closer to the overall median weight, 0.69). The LS estimator  remains different from the robust estimators, possibly due to the existence of casewise outliers in the data. MM-regression \citep{yohai:1985} behaves similarly to 2S-regression in this example. 
\begin{table}[!t]
\centering
\caption{Estimates and p-values of the regression coefficients for the imputed Boston Housing data.}\label{boston-estimates-imputedver}
\footnotesize
\begin{tabular}{lrrrcrrrr}
  \hline
Variable & \multicolumn{2}{c}{3S} & \multicolumn{2}{c}{ShootS} & \multicolumn{2}{c}{2S} & \multicolumn{2}{c}{LS}  \\ 
\cmidrule(lr){2-3}\cmidrule(lr){4-5}\cmidrule(lr){6-7}\cmidrule(lr){8-9}
 & Coeff. & P-Val. & Coeff. & P-Val. & Coeff. & P-Val. & Coeff. & P-Val. \\
  \hline
  log(lstat) & -0.243 & 0.000 & -0.264 & - & -0.227 & $<$0.001 & -0.385 & $<$0.001 \\ 
  rm$^2$ & 0.015 & 0.000 & 0.013 & - & 0.014 & $<$0.001 & 0.009 & $<$0.001 \\ 
  tax & -0.051 & 0.000 & -0.030 & - & -0.047 & $<$0.001 & -0.032 & 0.002 \\ 
  log(dis) & -0.125 & 0.000 & -0.161 & - & -0.129 & $<$0.001 & -0.144 & $<$0.001 \\ 
  ptratio & -0.026 & 0.000 & -0.028 & - & -0.025 & $<$0.001 & -0.027 & $<$0.001 \\ 
  nox$^2$ & -0.578 & 0.013 & -0.522 & - & -0.619 & 0.010 & -0.479 & $<$0.001 \\ 
  age & -0.023 & 0.645 & -0.037 & - & -0.037 & 0.471 & 0.051 & 0.386 \\ 
  black & -0.726 & 0.398 & 0.371 & - & -0.882 & 0.376 & -0.206 & 0.519 \\ 
  log(crim) & -0.006 & 0.513 & -0.001 & - & -0.012 & 0.233 & -0.012 & 0.213 \\ 
   \hline
\end{tabular}
\end{table}
\begin{table}[!t]
\center
\small
\caption{Pairwise squared norm distances between the estimates for the imputed Boston housing data.}\label{tab:boston-pairwise-AND-experiment}
\begin{tabular}{lcccc}
\hline
 & 3S & ShootS & 2S & LS \\ 
\hline
3S & -  & 0.862 & 0.172 & 5.486 \\
ShootS &  & - & 1.158 & 3.992 \\
2S &  &  & - & 6.366 \\
LS &  &  &  & - \\
\hline
\end{tabular}
\end{table}

\newpage
\section*{References}

\bibliographystyle{elsarticle-harv} 
\bibliography{3Sregression}